\documentclass[11pt,reqno]{amsart}

\theoremstyle{plain}
\newtheorem{thm}{Theorem}[section]
\newtheorem{lem}[thm]{Lemma}
\newtheorem{cor}[thm]{Corollary}
\newtheorem{prop}[thm]{Proposition}
\newtheorem{fact}[thm]{Fact}
\newtheorem{q}{Question}

\theoremstyle{definition}
\newtheorem{defn}[thm]{Definition}
\newtheorem{rem}[thm]{Remark}

\newcommand{\w}{\ensuremath{\omega}}
\newcommand{\om}{\w}
\renewcommand{\a}{\ensuremath{\alpha}}
\newcommand{\al}{\a}
\newcommand{\ap}{\a}
\renewcommand{\b}{\ensuremath{\beta}}
\newcommand{\be}{\b}
\newcommand{\bt}{\b}
\newcommand{\g}{\ensuremath{\gamma}}
\newcommand{\gm}{\g}
\newcommand{\e}{\ensuremath{\varepsilon}}
\newcommand{\vep}{\e}
\newcommand{\ee}{\e}
\newcommand{\ba}{\mbox{$\beta\!\cdot\!\alpha$}}
\newcommand{\ib}[1]{I_b(#1)}

\newcommand{\llt}{\mbox{$\ell_1$-tree}}
\newcommand{\llkt}{\mbox{$\ell_1$-$K$-tree}}
\newcommand{\llbbt}{\mbox{$\ell_1$-block} basis tree}

\newcommand{\llkbbt}{\mbox{$\ell_1$-$K$-block} basis tree}

\newcommand{\vsb}{\vspace{\baselineskip}}
\newcommand{\spp}{\operatorname{supp}}

\newcommand{\tsb}[1][\alpha]{\ensuremath{T({\mathcal S}_{#1},1/2)}}

\newcommand{\N}{\ensuremath{\mathbf{N}}}

\newcommand{\cf}{\ensuremath{\mathcal{F}}}

\newcommand{\cs}{\ensuremath{\mathcal{S}}}

\begin{document}

\title[Concerning the Bourgain $\ell_1$ index of a Banach
  space]{Concerning the Bourgain $\boldsymbol{\ell_1}$ index of a
  Banach space}  
\author{Robert Judd}
\address{Department of Mathematics\\
  University of Texas at Austin\\
  Austin, TX 78712-1082 \\
  U.S.A. }
\email{rjudd@math.utexas.edu}
\author{Edward Odell}
\address{Department of Mathematics
  University of Texas at Austin\\
  Austin, TX 78712-1082 \\
  U.S.A. }
\email{odell@math.utexas.edu}
\thanks{Research supported by the NSF and TARP}
\date{}

\begin{abstract}
  A well known argument of James yields that if a Banach space $X$
  contains $\ell_1^n$'s uniformly, then $X$ contains $\ell_1^n$'s
  almost isometrically.  In the first half of the paper we extend this
  idea to the ordinal $\ell_1$-indices of Bourgain.  In the second
  half we use our results to calculate the $\ell_1$-index of certain
  Banach spaces.  Furthermore we show that the $\ell_1$-index of a
  separable Banach space not containing $\ell_1$ must be of the form
  $\omega^{\alpha}$ for some countable ordinal \al.
\end{abstract}

\maketitle

\section{Introduction}

It is well known that if $\ell_p\ (1\leq p<\infty)$ or $c_0$ is
crudely finitely representable in a Banach space $X$, then it is
finitely representable in $X$.  This was shown for $\ell_1$ and $c_0$
by R.C. James \cite{j} and for $\ell_p\ (1<p<\infty)$ it is a
consequence of Krivine's theorem \cite{k} as noted by Rosenthal
\cite{r}, \cite{l}.  We may state this as

\begin{quotation}
  
  For all $p\in[1,\infty]$, every $K \geq1$, each $m\geq1$, and every
  $\vep>0$, there exists $n$ such that if $(x_i)_1^n$ is a normalized
  basic sequence in a Banach space $X$ with
  $(x_i)_1^n\stackrel{K}{\sim}\mbox{uvb }\ell_p^n$, then there exists
  a normalized block basis $(y_i)_1^m$ of $(x_i)_1^n$ satisfying
  $(y_i)_1^m\stackrel{1+\vep}{\sim}\mbox{uvb }\ell_p^m$.

\end{quotation}

Separable Banach spaces not containing $\ell_1$ may differ in the
complexity of $\ell_1^n$'s embedded inside.  This complexity is
measured in part by Bourgain's $\ell_1$-index \cite{b}.  Bourgain
considered trees $T(X,K)$ whose nodes are finite basic sequences in
the unit ball of a Banach space $X$, $K$-equivalent to the unit vector
basis of some finite dimensional $\ell_1$, for a fixed $K$.  The
$\ell_1$-$K$-ordinal index of $X,\ I(X,K)$, was then defined to be the
supremum of the orders of such trees.

The definition of the $\ell_1$-trees constructed by Bourgain may be
extended to $\ell_p$-trees $(1<p\leq\infty)$ (we explain all the
unfamiliar terms in the next section).  We extend the results on
finite representability of $\ell_p$ in $X$ to $\ell_p$-trees for $p=1$
or $\infty$.  We prove the following theorem in
Section~\ref{proof-main}.

\begin{thm}\label{theorem}
  
  For $p=1\mbox{ or }\infty$, for each $K>1$, for every
  $\alpha<\omega_1$, and any $\varepsilon>0$, there exists
  $\beta<\omega_1$ such that for all Banach spaces $X$, if $T$ is an
  $\ell_p$-tree on $X$ with constant $K$ and order, $o(T)\geq\beta$,
  then there exists an $\ell_p$ block subtree $T'$ of $T$ with
  constant $1+\varepsilon$ and order, $o(T')\geq\alpha$.\end{thm}

This theorem is not true in general for $1<p<\infty$, and in the final
section we explain why not.  We also show how the same ideas may be
applied to the $\ell_1$-${\mathcal S}_{\alpha}$-spreading models
introduced by Kiriakouli and Negrepontis \cite{kn}.

In Section~\ref{index} we apply our results to the problem of
calculating Bourgain's $\ell_1$-index $I(X)$ of certain spaces $X$.
We show for example that if $X$ is Tsirelson's space, then
$I(X)=\w^{\w}$.  We prove that $I(X)$ is always of the form $\w^{\ap}$
and relate $I(X)$ to the ``block'' Bourgain $\ell_1$-index $\ib{X}$
for spaces with a basis.  Both indices are defined in
Section~\ref{index}.

\section{Preliminaries on trees}

By a {\em tree}\/ we shall mean a countable, non-empty, partially
ordered set $(T,\leq)$ for which the set $\{ y\in T : y<x\}$ is
linearly ordered and finite for each $x\in T$.  The elements of $T$
are called {\em nodes}.  The {\em predecessor node}\/ of $x$ is the
maximal element $x'$ of the set $\{ y\in T : y<x\}$, so that if $y<x$,
then $y\leq x'$.  The {\em initial nodes}\/ of $T$ are the minimal
elements of $T$ and the {\em terminal nodes}\/ are the maximal
elements.  A {\em subtree}\/ of a tree $T$ is a subset of $T$ with the
induced ordering from $T$.  This is clearly again a tree.  Further, if
$T'\subset T$ is a subtree of $T$ and $x\in T$, then we write $x<T'$
to mean $x<y$ for every $y\in T'$. We will also consider trees related
to some fixed set $X$.  \emph{A tree on a set $X$} is a subset
$T\subseteq\cup_{n=1}^{\infty}X^n$ with the ordering given by:
$(x_1,\dots,x_m)\leq(y_1,\dots,y_n)$ if $m\leq n$ and $x_i=y_i$ for
$i=1,\dots,m$.

The property of trees which is most interesting here is their {\em
  order}.  Before we can define this we must recall some terminology.
Let the {\em derived tree}\/ of a tree$T$ be $D(T)=\{ x\in T : x<y
\mbox{ for some } y\in T\}$.  It is easy to see that this is simply
$T$ with all of its terminal nodes removed.  We then associate a new
tree $T^{\alpha}$ to each ordinal $\alpha$ inductively as follows.
Let $T^0=T$, then given $T^{\alpha}$ let $T^{\alpha+1}=D(T^{\alpha})$.
If $\alpha$ is a limit ordinal, and we have defined $T^{\beta}$ for
all $\beta<\alpha$, let $T^{\alpha}=\cap_{\beta<\alpha}T^{\beta}$.  A
tree $T$ is \emph{well-founded} provided there exists no subset
$S\subseteq T$ with $S$ linearly ordered and infinite.  The order of a
well-founded tree $T$ is defined as $o(T)=\inf\{\alpha :
T^{\alpha}=\emptyset \}$.

A tree $T$ on a topological space $X$ is said to be \emph{closed}
provided the set $T\cap X^n$ is closed in $X^n$, endowed with the
product topology, for each $n\geq1$.  We have the following result
(see~\cite{b},~\cite{d}) concerning the order of a closed tree on a
Polish space.
\begin{prop} \label{prop:wf=>cble} 
  If $T$ is a well-founded, closed tree on a Polish 
  (separable, complete, metrizable) space, then $o(T)<\omega_1$.
\end{prop}

A map $f:T\rightarrow T'$ between trees $T$ and $T'$ is a {\em tree
  isomorphism}\/ if $f$ is one to one, onto and order preserving.  We
will write $T\simeq T'$ if $T$ is tree isomorphic to $T'$ and
$f:T\stackrel{\sim}{\rightarrow}T'$ to denote an isomorphism.  From
now on we shall simply write {\em isomorphism}\/ rather than tree
isomorphism.

\begin{defn}  \label{min-tree-def}
  Minimal tree\end{defn} A tree $t$ is a {\em minimal tree of order
  \al}, for some ordinal $\al<\omega_1$, if for each tree $T$ of order
\al\ there exists a subtree $T'\subset T$ of order \al\ which is
isomorphic to $t$.  Notice that if $t$ is a minimal tree of order
\al\, then any subtree of $t$ of order \al\ is also a minimal tree of
order \al.  We construct certain minimal trees for each ordinal
$\alpha<\omega_1$ in Section~\ref{ord-trees}. \vspace{\baselineskip}

If $X$ is a Banach space and $(x_i)_1^m\subset X$ with $\|x_i\|=1\ 
(i=1,\dots,m)$ we write $(x_i)_1^m\stackrel{K}{\sim}\mbox{uvb }\ell_p^m$
if there exist constants $c,C$ with $c^{-1}C\leq K$ and
\[ c\Bigl(\sum_1^m|a_i|^p\Bigr)^{\frac{1}{p}} \leq
\Bigl\|\sum_1^ma_ix_i\Bigr\| \leq
C\Bigl(\sum_1^m|a_i|^p\Bigr)^{\frac{1}{p}} \] 
for all $(a_i)_1^m\subset {\mathbf R}$.

\begin{defn} \label{lpk-tree}
$\ell_p$-$K$-tree\end{defn}

An \emph{$\ell_p$-$K$-tree} on a Banach space $X$ is a tree $T$ on $X$
such that $T\subseteq\cup_{n=1}^{\infty}S(X)^n$ and
$(x_i)_1^m\stackrel{K}{\sim} \mbox{uvb } \ell_p^m$ for each
$(x_1,\dots,x_m)\in T$.  We say $T$ is an \emph{$\ell_p$-tree} on $X$
if $T$ is an $\ell_p$-$K$-tree for some $K$.  For $p=1$ this
definition is slightly different to that in \cite{b} where an
$\ell_1$-$K$-tree is the largest tree of this form.  In fact our trees
are subtrees of those.

\begin{defn}  \label{b-sub-tree}
Block subtree\end{defn} 

Let $T$ be an \llt\ on a Banach space $X$.  We say $S$ is a {\em block
  subtree}\/ of $T$, written $S\preceq T$, if $S$ is a tree on $X$
such that there exists a subtree $T'\subset T$ and an isomorphism $f:
T'\stackrel{\sim}{\rightarrow}S$ satisfying:
\begin{itemize}
  \item For each $x=(x_1,\dots,x_n)\in T'$, let $y=(x_1,\dots,x_m)$ be
    the initial node of $T'$ with $y\leq x$.  If $y$ is also an
    initial node of $T$, then let $k=0$, otherwise let
    $(x_1,\dots,x_k)$ be the predecessor node of $y$ in $T$.  Then
    $f(x)$ is a normalized block basis of $(x_{k+1},\dots,x_n)$.
  \item If $x=(x_1,\dots,x_n)\in T'$ has predecessor node
    $(x_1,\dots,x_m)$ in $T'$, then $f(x)=f(x')\cup(y_i)_1^k$, where
    $(y_i)_1^k$ is a normalized block basis of $(x_{m+1},\dots,x_n)$.
\end{itemize}
For each node $(y_1,\dots,y_k)=f(x)\in S$ we call $x$ the {\em parent
  node}\/ of $(y_1,\dots,y_k)$.  Note that if $T$ is an
$\ell_1$-$K$-tree on $X$ and $S$ is a block subtree of $T$ then $S$ is
also an $\ell_1$-$K$-tree on $X$.  However, if $T$ is an
$\ell_p$-$K$-tree on $X$ for $p>1$ and $S$ is a block subtree of $T$,
then $S$ is an $\ell_p$-$K^2$-tree on $X$.

\section{Ordinal Trees}\label{ord-trees}

Most of the work needed to prove Theorem~\ref{theorem} is concerned
with constructing certain general trees consisting of collections of
finite subsets of ordinals ordered by inclusion.  We first construct
specific minimal trees $T_{\alpha}$ of order \al\ for every ordinal
$\alpha<\omega_1$.  Once this is done we construct ``replacement
trees'' $T(\alpha,\beta)$ which are formed by replacing each node of
$T_{\alpha}$ by one or more copies of $T_{\beta}$, and show that
$T(\alpha,\beta)$ is a minimal tree of order \ba.  This gives us in
some sense a ``tree within a tree'' or ``an $\alpha$ tree of $\beta$
trees''.

These two results are used as follows: Given an arbitrary
$\ell_1$-$K$-tree on $X$ with $o(T)\geq\alpha^2$ we can find a subtree
isomorphic to $T(\alpha,\alpha)$.  For one of the $\alpha$ trees
inside this we either have a good constant---in which case we are
finished---or we take a vector in the linear span of one of its nodes
with a bad constant.  Putting some of these vectors together yields a
block subtree of order \a, each of whose nodes is ``bad'', and then
following the original argument of James these vectors together have a
good constant.

We now define the trees $T_{\alpha}$, and prove in
Lemma~\ref{t-alpha-min} that $T_{\alpha}$ is minimal of order.

\begin{defn}  \label{t-alpha-def}
  Minimal trees, $T_{\alpha}$
\end{defn} 
We define trees $T_{\alpha}$ of order $\alpha$ for each countable
ordinal $\alpha$ as subsets of $[1,\gamma]^{<\omega}$ ordered by
inclusion, for some ordinal $\gamma=\gamma(\alpha)<\omega_1$ where, if
$S$ is any set, then $[S]^{<\omega}$ is the collection of all finite
subsets of $S$.  We choose $\gamma(\alpha)$ and $T_{\alpha}$ by
induction as follows: Let $T_1=\{\{1\}\}$.  Given
$T_{\alpha}\subset[1,\gamma]^{<\omega}$ for some ordinal
$\gamma<\omega_1$, let $T_{\alpha+1}=\{ A\cup\{\gamma+1\} : A\in
T_{\alpha}\}\cup\{\{\gamma+1\}\}$.  Note that for $\beta<\al,\ 
(T_{\alpha+1})^{\beta}=\{ A\cup\{\gamma+1\} : A\in
(T_{\alpha})^{\beta}\}\cup\{\{\gamma+1\}\}\mbox{ and }
(T_{\alpha+1})^{\alpha}=\{\{\gamma+1\}\}$.  Thus
$o(T_{\alpha+1})=\al+1$ as required.


Finally, to define $T_{\alpha}$ for $\alpha$ a limit ordinal, let
$\alpha_n\nearrow\alpha$ be a sequence of ordinals increasing to \a,
and let $T_{\alpha_n}\subset[1,\beta_n]^{<\omega}$ for some
$\beta_n<\omega_1$.  Let $\beta=\sup_n\beta_n$ and $\gamma_n=\beta+n$
for each $n$.  Let $\tilde{T}_{\alpha_n}=\{ A\cup\{\gamma_n\} : A\in
T_{\alpha_n}\}$ and let
$T_{\alpha}=\cup_1^{\infty}\tilde{T}_{\alpha_n}$, ordered by
inclusion.  Notice that $\tilde{T}_{\alpha_n}$ is the same tree as
$T_{\alpha_n}$ with the same order and structure, but the nodes have
simply been relabeled.  The reason for doing this is that nodes from
different trees are now incomparable, and so the union
$\cup_1^{\infty}\tilde{T}_{\alpha_n}$ is a disjoint union.

To give an idea of what these trees look like we will construct the
trees $T_n$ and $T_{\omega}$ explicitly.
\begin{align*}
  T_1 & =  \{\{1\}\} \\
  T_2 & =  \{\{1,2\},\{2\}\} \\
  T_3 & =  \{\{1,2,3\},\{2,3\},\{3\}\} \\
      & \ \ \vdots \\
  T_n & = \{\{1,2,3,\dots,n\},\{2,3,\dots,n\},
          \dots,\{n-1,n\},\{n\}\}\ . \\
  \intertext{Then to construct $T_{\omega}$ we use the trees
    $\tilde{T}_n\ (n\geq1)$ as described above.}
  \tilde{T}_1 & =  \{\{1,\omega+1\}\}  \\
  \tilde{T}_2 & =  \{\{1,2,\omega+2\},\{2,\omega+2\}\} \\
              & \ \ \vdots \\
  \tilde{T}_n & =  \{\{1,\dots,n,\omega+n\},\dots,\{n,\omega+n\}\} \\
  T_{\omega}  & =
              \{\{1,\omega+1\},\{1,2,\omega+2\},\{2,\omega+2\},\dots,
              \{1,\dots,n,\omega+n\},\dots,\{n,\omega+n\},\dots\} \ .
\end{align*}

\begin{lem} \label{limit-trees}
  Let $\alpha<\omega_1$ be a limit ordinal and $T$ be a countable tree
  of order $\alpha$.  Then there exist a sequence $(\alpha_n)$ of
  successor ordinals and a sequence $(t_n)$ of subtrees $t_n\subset T$
  with $\alpha=\sup_n\alpha_n,\ o(t_n)=\alpha_n\mbox{ and }
  T=\cup_1^{\infty}t_n$. Moreover the trees $(t_n)$ are mutually
  incomparable, ie.\ if $x\in t_n$ and $y\in t_m$ with $n\neq m$, then
  $x$ and $y$ are incomparable.\end{lem}
\begin{proof}
  Suppose that $T$ has only finitely many initial nodes; let these be
  $x_1,\dots,x_n$, and let $t_i=\{y\in T : y\geq x_i\}$.  Then
  $\al=o(T)=\max_{1\leq i\leq n}o(t_i)=o(t_{i_0})$ for some $i_0\leq
  n$.  Let $t=\{y\in t_{i_0} : y>x_{i_0} \}$ and let $\beta=o(t)$.
  Since $\{x_{i_0}\}$ is the unique initial node of $t_{i_0}$, it
  follows that $t_{i_0}=t\cup\{x_{i_0}\}$ and hence
  $(t_{i_0})^{\beta}=\{x_{i_0}\}$.  Thus $\al=o(t_{i_0})=\beta+1$, a
  successor ordinal, contradicting the assumption that \al\ is a limit
  ordinal.
  
  Thus $T$ must have infinitely many initial nodes; let these be
  $(x_n)_1^{\infty}$ and let $t_n=\{y\in T : y\geq x_n \},\ 
  \alpha_n=o(t_n)$.  Note that these trees are mutually incomparable
  since the nodes $(x_n)_1^{\infty}$ are incomparable. We find that
  $\alpha_n$ is a successor ordinal using the same argument as above
  and from the definition of the order of a tree we have that
  $o(T)=\sup_n o(t_n)$ and hence $\alpha=\sup_n\alpha_n$.
\end{proof}

\begin{lem} \label{t-alpha-min}
$T_{\alpha}$ is a minimal tree of order \a.\end{lem}
\begin{proof}  
  The order of $T_{\alpha}$ is clear from the construction; we prove
  here that if $T$ is any tree of order $\alpha<\omega_1$, then there
  exists a subtree $T'\subset T$ such that $T'$ is isomorphic to
  $T_{\alpha}$.  We use induction on $\alpha$, the order of $T$.  The
  result is obvious for $\alpha=1$.
  
  Suppose the lemma is true for the ordinal $\alpha<\omega_1$.  Let
  $T$ have order $\alpha+1$, and hence $T^{\alpha}\neq\emptyset$.  Let
  $x$ be a terminal node of $T^{\alpha}$ and let $\tilde{T}=\{y\in T :
  y>x\}$; then $o(\tilde{T})=\alpha$. By assumption there exists a
  subtree $\tilde{T'}$ of $\tilde{T}$ and an isomorphism
  $f:T_{\alpha}\stackrel{\sim}{\rightarrow}\tilde{T'}$.  Clearly
  $T'=\tilde{T'}\cup{\{x\}}$ is a subtree of $T$ of order $\alpha+1$
  and we can extend $f$ to
  $F:T_{\alpha+1}\stackrel{\sim}{\rightarrow}T'$ to show that $T'$ is
  isomorphic to $T_{\alpha+1}$ as follows.  Recall from
  Definition~\ref{t-alpha-def} that we obtain $T_{\alpha+1}$ from
  $T_{\alpha}$ by setting $T_{\alpha+1}=\{ a\cup\{\gamma+1\} : a\in
  T_{\alpha}\}\cup\{\{\gamma+1\}\}$.  Setting $F(\{\gamma+1\})=x$ and
  $F(a\cup\{\gamma+1\})=f(a)$ makes $F$ the required isomorphism.
  
  If $\alpha$ is a limit ordinal, let the lemma be true for all
  $\beta<\alpha$ and let $T$ have order $\alpha$.  By
  Lemma~\ref{limit-trees}, $T=\cup_1^{\infty}t_n$ where
  $o(t_n)=\beta_n$, $\alpha=\sup_n\beta_n$, each $\be_n$ is a
  successor ordinal, and the trees $(t_n)$ are mutually incomparable.
  Let $\alpha_n\nearrow\alpha$ be the sequence of ordinals increasing
  to $\alpha$, and let $\tilde{T}_{\alpha_n}$ be the trees, from the
  definition of the minimal tree $T_{\alpha}$,
  Definition~\ref{t-alpha-def}.  Let $(\beta_{r_n})$ be a subsequence
  of $(\beta_n)$ so that $\alpha_n\leq\beta_{r_n}$ for all $n$.  Each
  tree $t_{r_n}$ contains a subtree of order $\alpha_n$; hence, by
  assumption, for each $n$ there exists $t'_{r_n}\subset t_{r_n}$ and
  an isomorphism $f_n:
  T_{\alpha_n}\stackrel{\sim}{\rightarrow}t'_{r_n}$.  Using the
  notation of Definition~\ref{t-alpha-def} we define $\tilde{f}_n:
  \tilde{T}_{\alpha_n}\stackrel{\sim}{\rightarrow}t'_{r_n}$ by
  $\tilde{f}_n(a\cup\{\gamma_n\})=f_n(a)$. Let
  $T'=\cup_1^{\infty}t'_{r_n}$ and
  $f:T_{\alpha}\stackrel{\sim}{\rightarrow}T'$ be the function
  $f=\cup_1^{\infty}\tilde{f}_n$.
\end{proof}

\begin{rem}\label{rem-after-min-tree}
  It follows that if $T$ is a tree of order $\beta\geq\alpha$, then
  there exists a subtree $T'\subset T$ such that $T'$ is isomorphic to
  $T_{\alpha}$.\vspace{\baselineskip}
\end{rem}
We now construct the replacement trees $T(\alpha,\beta)$, for each
pair of ordinals \a, $\b<\w_1$, promised earlier.  First we construct
the trees by induction, then we prove that $T(\alpha,\beta)$ has order
\mbox{$\beta\!\cdot\!\alpha$}.  Finally we show that $T(\alpha,\beta)$
is isomorphic to a subtree of $T_{\beta\cdot\alpha}$ and hence is a
minimal tree of order \ba\ as required.  The key to all of these
proofs is to use induction on $\alpha$ for an arbitrary $\beta$.

\begin{defn} \label{rep-tree-def}
  Replacement trees
\end{defn} 
For each pair $\alpha,\beta<\omega_1$ we construct a tree
$T(\alpha,\beta)$ and a map $f_{\alpha,\beta}:
T(\alpha,\beta)\rightarrow T_{\alpha}$ satisfying:
\begin{enumerate}

\item[(i)] For each $x\in T_{\alpha}$ there exists $I=\{1\}$ or $\N$
  and trees $t(x,j)\simeq T_{\beta},\ j\in I$, such that 
  $f_{\alpha,\beta}^{-1}(x)=\cup_{j\in I}t(x,j)$ (incomparable union) with
  $I=\{1\}$ if \al\ is a successor ordinal and $x$ is the unique initial
  node, or $\be<\omega$, and $I=\N$ otherwise.

\item [(ii)] For each pair $a,b\in T(\alpha,\beta),\ a\leq b$ implies
$f_{\alpha,\beta}(a)\leq f_{\alpha,\beta}(b)$.

\end{enumerate}

For each $\be<\om_1$, let $T(1,\be)=T_{\beta}$ and
$f_{1,\beta}:T(1,\be)\rightarrow T_1$ be given by
$f_{1,\beta}(a)=\{1\}\ \forall a\in T(1,\be)$.  Let $\al<\om_1$ and
suppose we have defined $T(\alpha,\beta) \mbox{ and }
f_{\alpha,\beta}$ for each $\beta<\omega_1$.  Roughly speaking, what
we do to go from $\alpha$ to $\alpha+1$ is to take $T_{\beta}$ and
then after each of its terminal nodes we put a $T(\alpha,\beta)$ tree.
This will give us the required tree, but we have to ensure that it is
well defined and that we keep track of the order relation.

Recall from Definition~\ref{t-alpha-def} that $T_{\alpha+1}=\{
a\cup\{\gamma+1\} : a\in T_{\alpha}\}\cup\{\{\gamma+1\}\}$ for some
$\gamma<\omega_1$.  Let $\delta_1, \delta_2$ be countable ordinals
with $T(\alpha,\beta)\subset[1,\delta_1]^{<\omega},
T_{\beta}\subset[1,\delta_2]^{<\omega}$.  Define a map $\ \tilde{
  }:[1,\delta_2]\rightarrow[\delta_1+1,\delta_1+\delta_2]$ by
$\eta\mapsto\tilde{\eta}=\delta_1+\eta$.  For all ordinals
$\lambda,\mu,\nu$, we have
$\lambda+\mu=\lambda+\nu\Rightarrow\mu=\nu$, hence this map is one to
one.  Thus, if we define $\tilde{T}_{\beta}=\{\tilde{a}:a\in
T_{\beta}\}$, then $\tilde{T}_{\beta}\simeq T_{\beta}$ as the map $\ 
\tilde{ }\ $ is merely relabeling the nodes, but the trees 
$\tilde{T}_{\beta}\mbox{ and }T(\alpha,\beta)$ are now incomparable
since if $a\in \tilde{T}_{\beta},\ b\in T(\alpha,\beta)$, then $a\cap
b=\emptyset$.

Let $(\tilde{x}_n)_I\ (I=\{1\}$ or $\N)$ be the set of terminal nodes
of $\tilde{T}_{\beta}$, a sequence of incomparable nodes and let
\begin{gather*} 
  T(\alpha+1,\beta)=\cup_{n\in I}\{a\cup \tilde{x}_n : a\in
          T(\alpha,\beta)\}\cup\tilde{T}_{\beta} \\ 
  f_{\alpha+1,\beta}(x)= \left\{ 
    \begin{array}{ll}
          f_{\alpha,\beta}(a)\cup\{\gamma+1\} & x=a\cup\tilde{x}_n\ \
                (a\in T(\alpha,\beta)) \\ 
          \{\gamma+1\} & x\in \tilde{T}_{\beta} 
    \end{array} \right. \ .
\end{gather*}

We need to show that the map $f_{\alpha+1,\beta}$ satisfies the
required properties.  Let $y\in T_{\alpha+1}$.  If $y=\{\gamma+1\}$,
then $f_{\alpha+1,\beta}^{-1}(y)=\tilde{T}_{\beta}\simeq T_{\beta}$.
Otherwise $y=a\cup\{\gamma+1\}$ for some $a\in T_{\alpha}$ and hence
\begin{alignat*}{2}
  f_{\alpha+1,\beta}^{-1}(y) 
    &= \cup_{n\in I}\{b\cup\tilde{x}_n
         : b\in f_{\alpha,\beta}^{-1}(a)\} && \quad\\
    &= \cup_{n\in I}\cup_{i\in I'}t_{n,i} && \quad\text{where
         $t_{n,i}\simeq T_{\beta}$ and $I'=\{1\}$ or \N} \\
    &= \cup_{j\in I''}t(y,j) && \quad\text{where $t(y,j)\simeq
         T_{\beta}$ and $I''=\{1\}$ or \N}
\end{alignat*}
as required.  Furthermore, the $t(y,j)$'s are incomparable.  The
second property is clear.

If $\alpha$ is a limit ordinal, let $\alpha_n\nearrow\alpha$ be the
sequence of ordinals increasing to $\alpha$ from
Definition~\ref{t-alpha-def} and suppose we have constructed
$T(\alpha_n,\beta),\ f_{\alpha_n,\beta}$ for each $\alpha_n$.  Let
$T(\alpha_n,\beta)\subset [1,\delta_n]^{<\omega},\ 
\delta=\sup_n\delta_n<\omega_1,\ \mbox{ and set} \gamma_n=\delta+n$
for each $ n$.  Then, as in the definition of the minimal trees, let
$\tilde{T}(\alpha_n,\beta)=\{a\cup\{\gamma_n\}\ :\ a\in
T(\alpha_n,\beta)\},\ 
\tilde{f}_{\alpha_n,\beta}(a\cup\{\gamma_n\})=f_{\alpha_n,\beta}(a)$,
and let $T(\alpha,\beta)=\cup_1^{\infty}\tilde{T}(\alpha_n,\beta),\ 
f_{\alpha,\beta}=\cup_1^{\infty}\tilde{f}_{\alpha_n,\beta}$.

\begin{lem} \label{order-rep-trees}
  \mbox{$o(T(\alpha,\beta))=\beta\!\cdot\!\alpha$}.
\end{lem} 
\begin{proof}
  We proceed by induction on $\alpha$ for an arbitrary fixed $\beta$.
  The result is obvious for $\alpha=1$.
  
  Suppose $o(T(\alpha,\beta))=\mbox{$\beta\!\cdot\!\alpha$}$.  By the
  construction of $T(\alpha+1,\beta)$ we have that
  $(T(\alpha+1,\beta))^{\beta\cdot\alpha}=\tilde{T}_{\beta}$ and hence
  $o(T(\alpha+1,\beta))=
  \mbox{$\beta\!\cdot\!\alpha$}+\beta=\mbox{$\beta\!\cdot\!(\alpha+1)$}$.
  If $\alpha$ is a limit ordinal and
  $o(\tilde{T}(\alpha_n,\beta))=o(T(\alpha_n,\beta))=
  \mbox{$\beta\!\cdot\!\alpha_n$}$ for each $n$, where
  $T(\alpha,\beta)=\cup_1^{\infty}\tilde{T}(\alpha_n,\beta)$ from
  Definition~\ref{rep-tree-def}, then $o(T(\alpha,\beta))=\sup_n
  o(\tilde{T}(\alpha_n,\beta))=
  \sup_n\mbox{$\beta\!\cdot\!\alpha_n$}=\mbox{$\beta\!\cdot\!\alpha$}$.
\end{proof}

The last of our results on these specially defined trees is the
following:

\begin{lem} \label{rep-in-min}
  $T(\alpha,\beta)$ is a minimal tree of order \ba.
\end{lem} 
\begin{proof}
  Since \mbox{$o(T(\al,\be))=\be\!\cdot\!\al$} and
  $T_{\beta\cdot\alpha}$ is a minimal tree of order
  \mbox{$\be\!\cdot\!\al$}, then by Remark~\ref{rem-after-min-tree} it
  is sufficient to prove that $T(\alpha,\beta)$ is isomorphic to a
  subtree of $T_{\beta\cdot\alpha}$.  We prove this by induction on
  $\alpha$ for an arbitrary $\beta$.  The result is obvious for
  $\alpha=1$ since $T(1,\be)=T_{\beta}$.
  
  Suppose $T(\alpha,\beta)$ is isomorphic to a subtree of
  $T_{\beta\cdot\alpha}$ and hence is a minimal tree of order \ba.
  Now, \mbox{$o(T_{\beta\cdot(\alpha+1)})=\beta\!\cdot\!(\alpha+1)$}
  so $o((T_{\beta\cdot(\alpha+1)})^{\beta\cdot\alpha})=\beta$, thus
  since $T_{\beta}$ is minimal it is isomorphic to a subtree of
  $(T_{\beta\cdot(\alpha+1)})^{\beta\cdot\alpha}$.  But by
  construction $(T(\alpha+1,\beta))^{\beta\cdot\alpha}\simeq
  T_{\beta}$ and hence is isomorphic to a subtree $t_0$ of
  $(T_{\beta\cdot(\alpha+1)})^{\beta\cdot\alpha}$.  Let the
  isomorphism which sends $(T(\alpha+1,\beta))^{\beta\cdot\alpha}$
  onto $t_0\subseteq (T_{\beta\cdot(\alpha+1)})^{\beta\cdot\alpha}$ be
  $a\mapsto a'$, so that if $(x_n)_1^{\infty}$ are the terminal nodes
  of $(T(\alpha+1,\beta))^{\beta\cdot\alpha}$, then
  $(x'_n)_1^{\infty}$ are their images in
  $(T_{\beta\cdot(\alpha+1)})^{\beta\cdot\alpha}$, the terminal nodes
  of $t_0$, under this map.  Let $T(x'_n)=\{y\in
  T_{\beta\cdot(\alpha+1)} : y>x'_n\}\subset
  T_{\beta\cdot(\alpha+1)}$, then
  \mbox{$o(T(x'_n))\geq\beta\!\cdot\!\alpha$ for each $ n$}.  Now, by
  assumption, for each $n$ there exists a subtree $t_n$ of $T(x'_n)$
  isomorphic to $T(\alpha,\beta)$ and hence the subtree
  $\tilde{T}=(\cup_1^{\infty}t_n)\cup t_0$ of
  $T_{\beta\cdot(\alpha+1)}$ is isomorphic to $T(\alpha+1,\beta)$ as
  required.

  Let $\alpha$ be a limit ordinal with
  $T(\alpha,\beta)=\cup_1^{\infty}\tilde{T}(\alpha_n,\beta)$ via the
  construction in Definition~\ref{rep-tree-def}, and let
  $T(\alpha_n,\beta)$ be isomorphic to a subtree of
  $T_{\beta\cdot\alpha_n}$ for each $n$.  Then
  $\tilde{T}(\alpha_n,\beta)$ is isomorphic to a subtree of
  $\tilde{T}_{\beta\cdot\alpha_n}$ for all $n$, where
  $\tilde{T}_{\beta\cdot\alpha_n}=\{a\cup\{\gamma'_n\} : a\in
  T_{\beta\cdot\alpha_n}\}$ for some $\gamma'_n$, from
  Definition~\ref{t-alpha-def}, and hence $T(\alpha,\beta)$ is
  isomorphic to a subtree of
  $T_{\beta\cdot\alpha}=\cup_1^{\infty}\tilde{T}_{\beta\cdot\alpha_n}$
  as required.
\end{proof}

\section{Proof of Theorem~\ref{theorem}} \label{proof-main}

We have shown everything we need about trees on subsets of ordinals
and we now want to apply this to $\ell_1$-trees on a Banach space $X$.

\begin{defn} \label{block-def}
Block of an $\ell_1$-tree\end{defn}

Let $T'$ be a subtree of an $\ell_1$-tree $T$.  A {\em block}\/ of
$T'$ with respect to $T$ is a normalized vector $v$ in the linear span
of some node $x=(x_1,\dots,x_n)\in T'$ where either:
 \begin{itemize}
   \item  $x$ is an initial node of $T$,
   \item  the initial node of the subtree $\{y\in T' : y\leq x\}$ of
     $T$ is an initial node of $T$, or 
   \item  $(x_1,\dots,x_m)$ is the predecessor node in $T$ of the
     initial node of $\{y\in T' : y\leq x\}$ and $v$ is in the linear
     span of $(x_{m+1},\dots,x_n)$.
 \end{itemize}
 If $T'=T$, then a block of $T$ is simply any normalized vector $v$ in
 the linear span of any node $(x_{1},\dots,x_n)$ of $T$.

\begin{lem} \label{blocks}
  Let $T$ be a tree on $X$ of order \mbox{$\beta\!\cdot\!\alpha$}
  isomorphic to $T(\alpha,\beta)$, and let $F:T\rightarrow T_{\alpha}$
  be the map from Definition~\ref{rep-tree-def} satisfying, for all
  $x\in T_{\ap},\ F^{-1}(x)=\cup_IT_n(x)$, where $I=\{1\}$ or \N,
  $T_n(x)\simeq T_{\beta}$ and the $T_n(x)$'s are mutually
  incomparable.  For each $x\in T_{\alpha} \mbox{ and } n\geq1$, let
  $b(x,n)$ be a block of $T_n(x)$ with respect to $T$.  Then there
  exists a block subtree $T'$ of $T$ and an isomorphism
  $g:T'\stackrel{\sim}{\rightarrow}T_{\alpha}$ satisfying: for every
  pair $ a,b\in T_{\alpha}$, with $a<b$, there exist $ x_1<\dots<x_m
  \mbox{ in } T_{\alpha},\mbox{ integers } n_{x_1},\dots,n_{x_m}\ 
  \mbox{and}\ k<m\ \mbox{such that}\ g^{-1}(a)=(b(x_i,n_{x_i}))_1^k$
  and $ g^{-1}(b)=(b(x_i,n_{x_i}))_1^m$.
\end{lem}

This sounds very complicated but all it is saying is that if you have
a tree on $X$, isomorphic to a replacement tree $T(\alpha,\beta)$,
then you can replace each $\beta$-subtree by a normalized vector in
the linear span of a node of that tree, and refine to get a block
subtree of order $\alpha$.

\begin{proof}
  As usual we prove this by induction on $\alpha$ for an arbitrary
  $\beta$.  The result is obvious for $\alpha=1$ and the only
  non-obvious case is the successor case.
  
  Assume that the lemma is true for $\alpha$.  Let $T$ be a tree on
  $X$ of order \mbox{$\beta\!\cdot\!(\alpha+1)$} isomorphic to
  $T(\alpha+1,\beta)$, let $F:T\rightarrow T_{\alpha+1}$ be the map
  with $F^{-1}(x)=\cup_IT_n(x)$ where $T_n(x)\simeq T_{\beta}$, and
  let $b(x,n)$ be given for each $x\in T_{\alpha+1},\ n\in I$.
 
  By construction of the replacement trees,
  $T^{\beta\cdot\alpha}\simeq T_{\beta}$, and in fact
  $T^{\beta\cdot\alpha}=F^{-1}(\{\gamma+1\})=T_1(\{\gamma+1\})$, where
  $T_{\alpha+1}=\{ a\cup\{\gamma+1\} : a\in
  T_{\alpha}\}\cup\{\{\gamma+1\}\}$ from Definition~\ref{t-alpha-def}.
  After each terminal node of $T^{\beta\cdot\alpha}$ lies a tree
  isomorphic to $T(\alpha,\beta)$.  Let these trees be
  $(t_j)_1^{\infty}$.  Let $j_0\geq1$ be such that
  $t_{j_0}>b(\{\gamma+1\},1)$; then $t_{j_0}\simeq T(\alpha,\beta)$
  and so the lemma applies giving us a block subtree $t'_{j_0}\preceq
  t_{j_0}$ and $g:t'_{j_0}\stackrel{\sim}{\rightarrow}T_{\alpha}$ as
  in the statement.  Now let
  \[T'=\{ (b(\{\gamma+1\},1),u_1,\dots,u_m) : (u_i)_1^m\in
  t'_{j_0}\}\cup\{(b(\{\gamma+1\},1))\}\] and let
  $G:T'\stackrel{\sim}{\rightarrow}T_{\alpha+1}$ by
  \[ G(a)=  \left\{ \begin{array}{ll}
                g((u_i)_1^m)\cup\{\gamma+1\} &  
                     a=(b(\{\gamma+1\},1),u_1,\dots,u_m) \\
                \{\gamma+1\} & a=(b(\{\gamma+1\},1))
                                  \end{array}
                          \right. \]
  then $G,T'$ clearly satisfy the lemma.  
  
  The proof where $\alpha$ is a limit ordinal just involves taking the
  union of the previous trees and functions.
\end{proof}

\begin{defn} \label{restr-subt} Restricted subtree of a tree.\end{defn} 
Let $T$ be a tree on a set $X$ and let $T'$ be a subtree of $T$.  We
define another tree on $X$, the \emph{restricted subtree $R(T')$ of
  $T'$ with respect to $T$}.  Let $x=(x_i)_1^n\in T'$ and let $y$ be
the unique initial node of $T'$ such that $y\leq x$; let $m\leq n$ be
such that $y=(x_i)_1^m$.  If $y$ is also an initial node of $T$, then
set $k=0$, otherwise let $k<m$ be such that $(x_i)_1^k$ is the
predecessor node of $y$ in $T$.  Finally, setting
$R(x)=(x_{k+1},\dots,x_n)$, we define $R(T')=\{ R(x) : x\in T'\}$.  It
is easy to see that $R(T')$ is isomorphic to $T'$.

\begin{proof}[Proof of Theorem~\ref{theorem} for $p=1$]
  Let $T$ be an $\ell_1$-$K$-tree of order $\alpha^2$ on $X$.  We show
  that there exists $T'\preceq T$ such that $T'$ is an
  $\ell_1$-$\scriptstyle \sqrt{K}$-tree of order $\alpha$.
  
  By Lemmas~\ref{t-alpha-min} and~\ref{rep-in-min}, $T(\alpha,\alpha)$
  is isomorphic to a subtree of $T$ and so we may assume that in fact
  $T(\alpha,\alpha)$ is isomorphic to $T$.  Now let $F:T\rightarrow
  T_{\alpha}$ be the map from Definition~\ref{rep-tree-def} with
  $F^{-1}(x)=\cup_IT_n(x),\ T_n(x)\simeq T_{\alpha}$ for every $ x\in
  T_{\alpha}$ and $ n\in I$.
  
  For each $x\in T_{\alpha}$ and $ n\geq1$ let
  $\tilde{T}_n(x)=R(T_n(x))$, the restriction being with respect to
  $T$.  Note that $\tilde{T}_n(x)$ is an $\ell_1$-$K$-tree isomorphic
  to $T_n(x)$.  If there exist $x\in T_{\alpha},\ n\in I$ such that
  $\tilde{T}_n(x)$ is an $\ell_1$-$\scriptstyle \sqrt{K}$-tree we are
  finished, since $\tilde{T}_n(x)$ has order $\alpha$.  Otherwise let
  $(x_1,\dots,x_m)$ be a node of $\tilde{T}_n(x)$ which is not
  $\scriptstyle \sqrt{K}$ equivalent to the unit vector basis of
  $\ell_1^m$ and let $b(x,n)=\sum_1^ma_ix_i$ where $(a_i)_1^m\subset
  {\mathbf R},\ \sum_1^m|a_i|>{\scriptstyle \sqrt{K}}$ and
  $\|b(x,n)\|=1$.  Note that $\|b(x,n)\|$ is a block of $T_n(x)$ with
  respect to $T$.

  By Lemma~\ref{blocks} there exists $T'\preceq T$ of order $\alpha$
  whose nodes are $(b(x_i,n_{x_i}))_1^m$ for some $n_{x_i}$ where $\{
  x_1<\dots<x_m \}=\{ y\in T: y\leq x \}$ for each $x\in T_{\alpha}$.
  We need only show that this tree has constant $\scriptstyle
  \sqrt{K}$.  Let $(y_i)_1^n$ be a node in $T'$ with parent node
  $x=(x_1,\dots,x_m)\in T$. Thus there exist subsets $E_i\subset
  \{1,\dots,m\},\ E_1<\dots<E_n$ (where $E<F$ means $\max E<\min F$)
  such that $y_i=\sum_{k\in E_i}a_kx_k$ for each $i$ and satisfying:
  \[ 1=\|y_i\|=\Bigl\|\sum_{E_i}a_kx_k\Bigr\| <
  \frac{1}{\sqrt{K}}\sum_{E_i}|a_k| \ . \]
  
  Let $(b_i)_1^n\subset {\mathbf R}$, then
  \begin{align*}
    \Bigl\|\sum_{i=1}^nb_iy_i\Bigr\| &=  \Bigl \|\sum_{i=1}^nb_i\sum_{k\in
      E_i}a_kx_k \Bigr \|  \\
                &\geq \frac{1}{K}\sum_1^n|b_i|\sum_{k\in E_i}|a_k| \\
    &> \frac{1}{K}\sum_1^n|b_i|\!\cdot\!\sqrt{K}\\
                &= \frac{1}{\sqrt{K}}\sum_1^n|b_i|
  \end{align*}
  as required.  These last few lines are James' argument.
  
  Now, if we choose the smallest $n$ so that
  $K^{\frac{1}{2^n}}\leq1+\varepsilon$, then we can iterate this
  argument to show that if $T$ is an $\ell_1$-$K$-tree of order
  $\alpha^{2^n}$, then there exists $T'\preceq T$ such that $T'$ is an
  $\ell_1$-$(1+\varepsilon)$-tree of order $\alpha$, which proves the
  theorem.
\end{proof}

\begin{rem}\label{rems-after-thm}{\ }\vspace{-\baselineskip}\\
\begin{enumerate}
\item The proof of Theorem~\ref{theorem} for $p=\infty$ is very
  similar to that for $p=1$, except that given an
  $\ell_\infty$-$K$-tree $T$ on $X$ of order $\al^{2^n}$ for $n$
  sufficiently large, we choose a block subtree $T'\preceq T$ of order
  \al\ to obtain $\|\sum_1^na_ix_i\|\leq(1+\ee)\sup_i|a_i|$ for all
  nodes $(x_i)_1^n\in T'$, and then the lower estimate follows
  automatically according to \cite{j}.
  
\item The proof of the theorem also gives some fixed points---that is,
  ordinals $\alpha$ such that if we have an $\ell_1$-$K$-tree of order
  $\alpha$, then for any $\ee>0$ we can get a block subtree of this
  which is an $\ell_1$-$(1+\varepsilon) $-tree also of order $\alpha$.
  In fact we see from the proof that this is true for every countable
  ordinal $\alpha$ which satisfies $\be<\al$ implies $\be^n<\al$ for
  each $n\geq1$.  From basic results on ordinals we see that \al\ 
  satisfies this condition if and only if \al\ is of the form
  $\alpha=\om^{\omega^{\gamma}}$ for some ordinal $\gamma$ (see
  Fact~\ref{ords} below).
\end{enumerate}
\end{rem}

\section{Calculating the $\ell_1$ index of a Banach space} 
\label{index}

\begin{defn} Block basis tree\end{defn}
A {\em block basis tree}\/ on a Banach space $X$, with respect to a
basis $(e_i)_1^{\infty}$ for $X$, is a tree $T$ on $X$ such that every
node $(x_i)_1^n$ of $T$ is a block basis of $(e_i)_1^{\infty}$.
Moreover, if $T$ is also an $\ell_1$-$K$-tree, then we say $T$ is an
$\ell_1$-$K$-block basis tree.
\begin{defn} 
  The $\ell_1$-indices of a Banach space $X$: $I(X)$ and $I_b(X)$.
\end{defn} 
Let $X$ be a separable Banach space and for each $K\geq1$ set
\[I(X,K)=\sup\{ o(T) : T\mbox{ is an $\ell_1$-$K$-tree on }X\}\ .\] 
The \emph{Bourgain $\ell_1$-index} of $X$ \cite{b} is then given by
\[I(X)=\sup_{1\leq K<\infty}\{I(X,K)\}\ .\]  
By Bourgain, $I(X)<\w_1$ if and only if $X$ does not contain $\ell_1$.
  
The block basis index is the analogous index to $I(X)$ except that it
is only defined on block basis trees.  For a Banach space $X$ with a
basis $(e_i)$, and $K\geq1$, set
\[I_b(X,K,(e_i))=\sup\{ o(T) : T\mbox{ is an $\ell_1$-$K$-block
  basis tree w.r.t.\ $(e_i)$ on }X\}\ .\] 
The \emph{block basis index} is then given by
\[I_b(X,(e_i))=\sup\{I_b(X,K,(e_i)) : 1\leq K<\infty \}\ .\]
When the basis in question is fixed we shall write $I_b(X,K)$ rather
than $I_b(X,K,(e_i))$ etc.  It is worth recalling here that $I_b(X)$
is not in general independent of the basis.  It is clear, however,
that $I_b(X,K,(e_i))\leq I(X,K)$ for every $X$, $K$ and $(e_i)$.

We next state some facts about ordinals.  The proofs may be found in
Monk \cite{m}.

\begin{fact} \label{ords}
  Let \al\ be an infinite countable ordinal.  Then the following
  statements hold:
\begin{itemize}
\item[(i)] There exist $k\geq1$, (countable) ordinals
  $\theta_1>\dots>\theta_k\geq0$ and $n_i\geq1\ (i=1,\dots,k)$,
  uniquely determined by \al\ such that
  \mbox{$\al=\om^{\theta_1}\!\cdot\!
    n_1+\dots+\om^{\theta_k}\!\cdot\! n_k$}.  This is the \emph{Cantor
    normal form} of an ordinal.
  
\item[(ii)] For all $\be<\al,\ \be\cdot2<\al$ if and only if there
  exists $\gamma<\om_1$ such that $\al=\om^{\gamma}$.
  
\item[(iii)] For all $\be<\al,\ \be^2<\al$ if and only if there exists
  $\gamma<\om_1$ such that $\al=\om^{\omega^{\gamma}}$.
                
\item[(iv)] If \mbox{$\al=\om^{\theta_1}\!\cdot\!
    n_1+\dots+\om^{\theta_k}\!\cdot\! n_k$}, then
  \mbox{$\om\!\cdot\!\al=\al$} if and only if $\theta_k\geq\om$.
                
\item[(v)] If \mbox{$\al=\om^{\theta_1}\!\cdot\!
    n_1+\dots+\om^{\theta_k}\!\cdot\! n_k$}, then
  \mbox{$\al\!\cdot\!\om=\om^{\theta_1+1}$}.

\end{itemize}
\end{fact}

Our first result of this section is to show how we may refine
$\ell_1$-trees in a Banach space with a basis to get $\ell_1$-block
basis trees, and explain how this relates to the indices.  We then
show that both $I(X)$ and $I_b(X)$ are of the form $\om^{\alpha}$ for
some \al, and that if $\al\geq\om$ for either index, then the indices
are the same.  The block basis trees are much easier to work with, and
once we have the block index of a space we have a good idea what the
index is.  In the second part of this section we use this idea to
calculate the index of some Tsirelson type spaces.\vsb

\noindent\textbf{Notation}\\
\noindent For a Banach space $X$ let $B(X)=\{ x\in X : \|x\|\leq1 \}$ and 
$S(X)=\{ x\in X : \|x\|=1 \}$ denote the unit ball and unit sphere of
$X$ respectively.  If $(x_i)_{i\in I}\subset X$, where $I\subset\N$,
let $[x_i]_{i\in I}$ be the closed linear span of these vectors.

If $X$ is a Banach space with basis $(e_i)_1^{\infty}$ let
$E_n=[e_i]_1^n$, let $P_n:X\rightarrow E_n$ be the basis projection
onto $E_n$ given by $P_n(\sum a_ie_i)=\sum_1^na_ie_i$, and let
$X_n=[e_i]_{n+1}^{\infty}$.  Finally, we define the support of $x\in
X$ with respect to $(e_i)_1^{\infty}$ as $\spp(x)=\{ n\geq1 :
(P_n-P_{n-1})(x)\neq0 \}$.  Thus, if $x=\sum_Fa_ie_i$ with $a_i\neq0$
for $i\in F$, then $\spp(x)=F$.  If $x=(x_1,\dots,x_n)$ is a sequence
of vectors, then $\spp(x)=\cup_1^n\spp(x_i)$.  In the following $X$
will always denote a separable Banach space not containing $\ell_1$.

\begin{prop} \label{ind-to-block-ind}
  Let $X$ have a basis, then \mbox{$I(X,K)\geq\om\!\cdot\!\al$}
  implies that $I_b(X,K+\varepsilon)\geq\ap$ for every $\ee>0$.
\end{prop}

We first prove the following elementary lemma:

\begin{lem} \label{block}
  Let $X$ be a Banach space with basis $(e_i)_1^{\infty}$ and let $T$
  be an \llt\ of order \om\ on $X$, then for each $ n\geq 1$ there
  exists a block $x$ of $T$ with $P_nx=0$.
\end{lem}

\begin{proof}
  There exists $m>n$ such that the linear space spanned by
  $(y_i)_1^m\in T$ has dimension greater than $n$.  Thus the
  restriction of $P_n$ to $[y_i]_1^m$ is not one to one and hence
  there exists $x\in[y_i]_1^m$ with $\|x\|=1$ and $P_nx=0$.
\end{proof}

\begin{proof}[Proof of Proposition~\ref{ind-to-block-ind}]  
  If \mbox{$I(X,K)\geq\w\!\cdot\!\ap$}, then there exists an \llkt\ 
  $T$ on $X$ of order \mbox{$\w\!\cdot\!\ap$} and this in turn, by
  Lemma~\ref{rep-in-min}, has an $\ell_1$-$K$-subtree $T'$ isomorphic
  to $T(\ap,\w)$.  Thus we may assume that $T$ itself is isomorphic to
  $T(\ap,\w)$.  We prove the following statement:
  \begin{quotation}
    For all $\ap<\w_1$, each $ l\geq0$, and every $\ee>0$, if $T$ is
    an \llkt\ isomorphic to $T(\ap,\w)$, then there exists an
    $\ell_1$-$K$-block subtree $T'$of $T$ of order \al\ such that for
    any node $(y_i)_1^m\in T'$ there exists $l=k(1)<\dots<k(m+1)$ with
    $\|y_i-P_{k(i+1)}y_i\|<\ee$ and $P_{k(i)}y_i=0\ (i=1,\dots,m)$.
  \end{quotation}
  
  We induct on \al; the statement is clear for $\al=1$ by
  Lemma~\ref{block}.  Suppose we have proved the statement for \al,
  and let $T\simeq T(\al+1,\om)$.  Let $F:T\rightarrow T_{\ap+1}$ be
  the map $F^{-1}(x)=\cup_IT_n(x)\ (I=\{1\}\mbox{ or }\N)$, from
  Definition~\ref{rep-tree-def}, such that $T_n(x)\simeq T_{\w}$ for
  each $n$ and $x$, and the $T_n(x)$'s are mutually incomparable.  Let
  $z$ be the unique initial node of $T_{\ap+1}$.  By Lemma~\ref{block}
  we can find a block $b(1,z)$ of $T_1(z)$ such that $P_lb(1,z)=0$ and
  we can find $l'>l$ such that $\|b(1,z)-P_{l'}b(1,z)\|<\ee$.  Let
  $\tilde{T}$ be a subtree of $T$ isomorphic to $T(\ap,\w)$ with
  $b(1,z)<\tilde{T}$.  Applying the induction hypothesis to
  $\tilde{T}$ we obtain $\tilde{T}'\preceq\tilde{T}$ such that
  $P_{l'}y_i=0$ for every node $(y_i)_1^m\in \tilde{T}'$.  Let $T'=\{
  (b(1,z),y_1,\dots,y_m) : (y_i)_1^m\in\tilde{T}' \}\cup\{ (b(1,z))
  \}$.  Then $T'$ is the required block subtree.
  
  Now let \al\ be a limit ordinal and suppose we have proved the
  statement for each $\be<\al$.  Let $(\al_n)$ be the sequence of
  ordinals increasing to \al\ such that $T=\cup_1^{\infty} T(n)$ where
  the trees $T(n)$ are mutually incomparable and $T(n)\simeq
  T(\al_n,\om)$.  Applying the hypothesis to each $T(n)$ we obtain
  block subtrees $T(n)'\preceq T(n)$.  Then $T'=\cup_1^{\infty}T(n)'$
  is the required block subtree.
  
  Thus, if we have an \llkt\ $T$ of order \mbox{$\w\!\cdot\!\ap$} and
  $\ee'>0$, then let $T'$ be the $\ell_1$-$K$-block subtree of $T$
  from above.  For each node $(y_i)_1^m$ of $T'$ let
  $(k(i))_1^{m+1}\subset{\mathbf N}$ be the sequence from above and
  let
  \[v_i=\frac{P_{k(i+1)}y_i}{\|P_{k(i+1)}y_i\|}\ (i=1,\dots,m)\ .\]
  The sequence $(v_i)_1^m$ is a uniform perturbation of a basis $K$
  equivalent to the unit vector basis of $\ell_1^m$.  Hence, if $\ee'$
  is chosen sufficiently small, then $(v_i)_1^m$ is $K+\ee$ equivalent
  to the unit vector basis of $\ell_1^m$.  This completes the proof
  since if we replace the nodes $(y_i)_1^m$ with $(v_i)_1^m$ as above,
  then we obtain $\tilde{T}$, a block basis tree of order \al\ and
  constant $(K+\ee)$, so that $\ib{X,K+\ee}\geq\ap$ as required.
\end{proof}

\begin{thm} \label{ibx=wa}
  Let $X$ be a Banach space with a basis, then $\ib{X}=\w^{\ap}$ for
  some $\ap<\w_1$.
\end{thm}
To prove this theorem we need some preliminary results.  We first show
that there is no \llkbbt\ whose order is the same as the block basis
index, and hence $I_b(X)$ must be a limit ordinal.  Then we show that
$\b<I_b(X)$ implies that \mbox{$\b\cdot2<I_b(X)$}, which completes the
proof.

\begin{lem} \label{ibx-lim}
  Let $X$ be a Banach space with a basis and $K\geq1$, then
  $\ib{X,K}\neq\ib{X}$.  In particular $\ib{X}$ is a limit ordinal.
\end{lem}

\begin{proof}
  We prove by induction on \al\ that for every Banach space $X$ with a
  basis and any $K\geq1$, if $\ib{X,K}=\ap$, then $\ib{X}>\ap$.  This
  is trivial for $\ap=1$.
  
  Let the result be true for \al\ and suppose, if possible, that it is
  false for $\ap+1$.  Let $X$ be a Banach space with basis
  $(e_i)_1^{\infty}$ and $K\geq1$ such that $\ib{X,K}=\ib{X}=\ap+1$.
  Now there exists an \llkbbt\ $T$ of order $\ap+1$ isomorphic to the
  minimal tree $T_{\ap+1}$.  Let $x=(x_1)$ be the unique initial node
  of $T$, let $k=\max(\spp x_1)$, let $X_k$ be the subspace of $X$
  spanned by $(e_i)_{i>k}$ and let $T(\ap)=\{ (y_i)_1^m :
  y=(x_1,y_1,\dots,y_m)\in T \mbox{ and } y>x \}$.  Clearly $T(\ap)$
  is an \llkbbt\ on $X_k$ of order \al, and so $\ib{X_k}>\ap$,
  otherwise $\ib{X_k,K}=\ap=\ib{X_k}$ contradicting our assumption.
  Thus there exists an \llbbt\ $T'$ on $X_k$ of order $\ap+1$ for some
  constant $K'\geq1$. But now the tree $\tilde{T}=\{
  (x_1,u_1,\dots,u_l) : (u_1,\dots,u_l)\in T' \}\cup\{(x_1)\}$ is an
  \llbbt\ on $X$ of order $\ap+2$ for some constant $K''$
  contradicting the assumption that $\ib{X}=\ap+1$.  This proves the
  result for $\al+1$.
  
  Let \al\ be a limit ordinal and suppose the result is true for every
  $\ap'<\ap$, but false for \al.  Again let $X$ be a Banach space with
  basis $(e_i)_1^{\infty}$, $K\geq1$ such that $\ib{X,K}=\ib{X}=\ap$
  and $T$ an \llkbbt\ of order $\ap$ isomorphic to the minimal tree
  $T_{\ap}$.  By Lemma~\ref{limit-trees} there exists a sequence of
  ordinals $(\ap_n)$ such that $\ap=\sup_n(\ap_n+1)=\sup_n\ap_n$ and
  mutually incomparable trees $t_n$ for each $n$ such that $t_n\simeq
  T_{\ap_n+1}$ and $T=\cup_nt_n$.  For each $n$ let
  $z_n=(w_i)_1^{k_n}$ be the unique initial node of $t_n$ and let
  $t'_n=\{ (y_i)_1^m : y=(w_1,\dots,w_{k_n},y_1,\dots,y_m)\in t_n
  \mbox{ and } y>z_n \}$, a tree isomorphic to $T_{\ap_n}$.  Clearly
  $T'=\cup_nt'_n$ is a tree on $X_1$ with order \al.  Let
  $\tilde{T}=\{ (e_1,u_1,\dots,u_l) : (u_1,\dots,u_l)\in
  T'\}\cup\{(e_1)\}$.  This is an \llbbt\ of order $\al+1$,
  contradicting the assumption that $\ib{X}=\ap$.  This proves the
  first part of the lemma.
  
  Suppose, if possible, that $\ib{X}=\al+1$ for some \al.  Then there
  exists an \llkbbt\ $T$ of order $\ap+1$ for some $K$ contradicting
  the previous result.
\end{proof}

\begin{lem} \label{tree-on-xn}
  Let $X$ be a Banach space with basis $(e_i)_1^{\infty}$.  If
  $\be<\ib{X}$, then there exists $K>1$ such that $\ib{X_n,K}\geq\be$
  for every $n\geq1$.
\end{lem}

\begin{proof}
  The result is trivial for $\be<\w$.  Suppose first that $\be<\ib{X}$
  is a limit ordinal and let $T$ be an \llkbbt\ on $X$ of order \be.
  Let $T(n)=\{(x_i)_{n+1}^l : \exists (x_i)_1^l\in T\mbox{ with
    }l>n\}$.  $T(n)$ is clearly a block subtree of $T$ and an \llkbbt\ 
  on $X_n$.  Moreover, $o(T(n))=\be$, otherwise $o(T)\leq
  o(T(n))+n<\be$, a contradiction.

  Now let $\be<\ib{X}$ be a successor ordinal greater than $\w$, then
  $\be=\be'+k$ for some limit ordinal $\be'\geq\w$ and $k\geq1$.  From
  the limit ordinal case there exists $K>1$ such that
  $\ib{X_m,K}\geq\be'$ for every $m\geq1$.  Now, $X$ contains
  $\ell_1^n$'s uniformly so there exists $m>n$ and a normalized block
  basis $(x_i)_1^k$ of $[e_i]_{n+1}^m$ which is $K$ equivalent to the
  unit vector basis of $\ell_1^k$.  Let $T$ be an \llkbbt\ on $X_m$ of
  order $\be'$ and let $T(n)=\{ (x_1,\dots,x_k,u_1,\dots,u_l) :
  (u_1,\dots,u_l)\in T\}\cup\{ (x_1,\dots,x_k),\dots,(x_1) \}$.  Then
  $T(n)$ is an \llbbt\ on $X_n$ of order $\bt'+k=\be$ and some
  constant which depends only on $K$.
\end{proof}

\begin{proof}[Proof of Theorem~\ref{ibx=wa}]
  We show that if $\be<\ib{X}$, then $\be\cdot2<\ib{X}$, which is
  enough to prove the theorem by Fact~\ref{ords} (ii).  Let
  $\be<\ib{X}$ and $T$ be an \llkbbt\ on $X$ of order \be.  For each
  $n$ let $T(n)$ be an \llkbbt\ on $X_n$ of order \be\ from
  Lemma~\ref{tree-on-xn}.  Let $(a_i)$ be the collection of terminal
  nodes of $T$ and for each $i\geq1$ let $n(i)=\max(\spp a_i)$.
  Finally, setting $\tilde{T}(n(i))=\{ a_i\cup x : x\in T(n(i)) \}$,
  we have that $\tilde{T}=T\cup(\cup_i\tilde{T}(n(i)))$ is an \llbbt\ 
  of order \mbox{$\be\!\cdot\!2$} and hence
  \mbox{$\ib{X}>\be\!\cdot\!2$} as required.
\end{proof}

\begin{thm} \label{ix=wa}
  Let $X$ be a separable Banach space, then $I(X)=\w^{\ap}$ for some
  $\ap<\w_1$.
\end{thm}
The proof of this theorem is similar to that of Theorem~\ref{ibx=wa},
but without a basis for $X$ we have to work harder.

\begin{lem} \label{max-nodes}
  Let $T$ be a countable tree of order $\al<\w_1$, $M$ the collection
  of maximal nodes of $T$, $M=\cup_{i=1}^nM_i$ a partition of $M$, and
  $T_i=\{ x\in T : x\leq m\mbox{ for some }m\in M_i\}$. Then
  $o(T_i)=\al$ for some $1\leq i\leq n$.
\end{lem}

\begin{proof}
  We prove by induction on \al.  The result is obvious for $\al=1$.
  Suppose it is true for \al, and let $T$ be a countable tree of order
  $\al+1$ with $M,\ M_i,\ T_i$ as above.  Let $(a_j)$ be the sequence
  of initial nodes of $T$ and $t_j=\{ x\in T : x\geq a_j \}$.  Clearly
  the $t_j$'s are mutually incomparable and $T=\cup_jt_j$, hence
  $o(t_{j_0})=\al+1$ for some $j_0$.  Let $t'=\{ x\in T : x>a_{j_0}
  \}$, then $o(t')=\al$.  Now, $M=\cup_{i=1}^nM_i$ also partitions the
  terminal nodes of $t'$ and setting $t'_i=\{ x\in t' : x\leq m \mbox{
    for some }m\in M_i \}$ we have $o(t'_{i_0})=\al$ for some $i_0$ by
  assumption.  Now $\{a_{j_0}\}\cup t'_{i_0}$ is a tree of order
  $\al+1$ and $\{a_{j_0}\}\cup t'_{i_0}\subseteq T_{i_0}$.  Thus
  $o(T_{i_0})=\al+1$ as required.
  
  Let \al\ be a limit ordinal and suppose the result is true for each
  $\al'<\al$.  Write $T=\cup t_k$ as a union of mutually incomparable
  trees $t_k$ of order $\al_k$ where $\sup_k\al_k=\al$.  Given $M,\ 
  M_i,\ T_i$ as above let $t_{k,i}=\{x\in t_k : x\leq m\mbox{ for some
    }m\in M_i\}$ and let $i(k)\in\{1,\dots,n\}$ satisfy
  $o(t_{k,i(k)})=\al_k$ for each $k$, by assumption.  Let $N_i=\{
  k\geq1 : i(k)=i \}$, then $N_i$ must be infinite for some $i_0$, so
  let $N_{i_0}=(k_j)_1^{\infty}$.  Now for each $j$ we have
  $t_{k_j,i(k_j)}=t_{k_j,i_0}\subseteq T_{i_0}$ and the trees
  $t_{k_j,i_0}$ are mutually incomparable, thus
  $o(\cup_jt_{k_j,i_0})=\al$ which implies $o(T_{i_0})=\al$ as
  required.
\end{proof}

\begin{lem} \label{ix-lim}
  Let $X$ be a separable Banach space not containing $\ell_1$ and
  $K\geq1$, then $ I(X,K)\neq I(X)$.  In particular $I(X)$ is a limit
  ordinal.
\end{lem}

\begin{proof}
  Let $I(X,K)=\ap$ for some $\ap<\w_1$ and $T$ be an \llkt\ on $X$ of
  order \al.  Recall that a Banach space $X$ is ${\mathcal L}_1$-$K$
  if there exists a collection $(E_n)$ of finite dimensional subspaces
  of $X$ with $d(E_n,\ell_1^{{\rm dim} E_n})\leq K$ for every $ n$,
  and for each finite set $F\subset X$ and all $\ee>0$ there exists
  $n$ such that the distance from $x$ to $E_n$ is less than $\ee$ for
  all $x$ in $F$.  Also recall that every infinite dimensional
  ${\mathcal L}_1$ space contains $\ell_1$.  See~\cite{lt} for more
  information on ${\mathcal L}_1$ spaces.
  
  Now let $M$ be the set of maximal nodes of $T$.  Clearly this
  defines a collection of finite dimensional subspaces $[x_i]_1^n$
  such that $d([x_i]_1^n,\ell_1^n)\leq K$, where $(x_i)_1^n\in M$.
  Thus, since $X$ doesn't contain $\ell_1$, it is not a ${\mathcal
    L}_1$ space and hence there exist $F=\{z_1,\dots,z_r\}\subseteq
  S(X)$ and $\ee>0$ such that for each $m=(x_i)_1^n\in M$ there exists
  $i(m)\in\{1,\dots,r\}$ with $d(z_{i(m)},S([x_i]_1^n))>\ee$.  For
  $i=1,\dots,r$ set $M_i=\{m\in M : i(m)=i\}$.  Then $M=\cup_1^rM_i$
  partitions $M$ and defines $T=\cup_1^rT_i$ as in
  Lemma~\ref{max-nodes}.  So, from the lemma, we have $o(T_{i_0})=\al$
  for some $i_0\leq r$.  Let $T'=\{ (z_{i_0},u_1,\dots,u_m) :
  (u_1,\dots,u_m)\in T_{i_0}\}\cup\{(z_{i_0})\}$, then this is an
  \llt\ on $X$, for some constant $K'=K'(K,\ee)$, of order $\al+1$.
  Thus $I(X)>\al=I(X,K)$ which completes the first part of the proof.
  The argument that $I(X)$ is a limit ordinal is the same as for
  $\ib{X}$.
\end{proof}

\begin{lem} \label{lem-for-ix=wa}
  Let $T$ be a tree on $X$ of order \al, where \al\ is a limit
  ordinal.  Let $F\subset S(X^{\ast})$ be finite and $X_F=\{ x\in X :
  x^{\ast}(x)=0\ \forall x^{\ast}\in F\}$.  Then there exists a block
  subtree $T'$ of $T$ with $o(T')=\al$ and $T'\subseteq X_F$.
\end{lem}

\begin{proof}
  Let $|F|=n$.  We note that \al\ is a limit ordinal if and only if
  \mbox{$\al=\w\!\cdot\!\be$} for some ordinal \be, and prove the
  lemma by induction on \be.
  
  For $\be=1,\ \al=\w$, and let $T$ be isomorphic to $T_{\w}$.  Notice
  that if $(x_i)_1^{n+1}\in T$, then there exists $x\in
  S([x_i]_1^{n+1})$ with $x\in X_F$.  Thus for each $k$ there exists a
  node $(x^k_i)_{i=1}^{l}\in T$, for $l$ sufficiently large, from
  which we may extract a normalized block basis $(y^k_j)_{j=1}^k$ of
  $(x^k_i)_{i=1}^{l}$ which is contained in $X_F$ and such that $T'=\{
  (y^k_1,\dots,y^k_j) : 1\leq j\leq k,\ k\geq1\}$ is a block subtree
  of $T$.  This is now the required tree.
  
  Suppose the result is true for $\be$ and let
  \mbox{$\al=\w\cdot(\be+1)=\w\!\cdot\!\be+\w$} and $T$ be a tree of
  order \al.  Since $T$ has a subtree isomorphic to $T_{\alpha}$ we
  may assume $T\simeq T_{\alpha}$.  Now $T^{\w\cdot\bt}$ is isomorphic
  to $T_{\w}$ and we apply the case $\be=1$ to obtain a block subtree
  $\tilde{T}\subset T^{\w\cdot\bt}$ of order $\w$, contained in $X_F$.
  Let $(\tilde{a}_i)$ be the sequence of terminal nodes in $\tilde{T}$
  and $a_i$ the parent node of $\tilde{a}_i$ in $T^{\w\cdot\bt}$ for
  each $i$.  Let $T(i)=\{ x\in T : x>a_i \}$, then
  \mbox{$o(T(i))\geq\w\!\cdot\!\be$}.  Thus we may apply the induction
  hypothesis to $R(T(i))$ (the restricted tree from
  Definition~\ref{restr-subt}) for each $i$ to obtain block subtrees
  $T(i)'\subset T(i)$ with \mbox{$o(T(i)')=\w\!\cdot\!\be$} and
  $T(i)'\subset X_F$. Finally, $T'=\tilde{T}\cup(\cup_iT(i)')$ is the
  required tree of order \al.
  
  Let \be\ be a limit ordinal and suppose the result is true for all
  $\be'<\be$.  Let $(\bt_n)$ be the increasing sequence of ordinals
  whose limit is \be, then
  \mbox{$\al=\w\!\cdot\!\be=\sup_n\w\!\cdot\!\be_n$} so that if $T$ is
  a tree of order \al\ and then $T$ contains mutually incomparable
  trees of order \mbox{$\w\!\cdot\!\be_n$} for each $\be_n$.  We apply
  the hypothesis to each of these trees to obtain the result.
\end{proof}

\begin{proof}[Proof of Theorem~\ref{ix=wa}]
  By Lemma~\ref{ix-lim}, $I(X)=\al$ for some limit ordinal \al.  We
  show that if $\be<\al$ is a limit ordinal, then
  \mbox{$\be\!\cdot\!2<\al$}.  It follows that if $\be<\al$ is a
  successor ordinal, then \mbox{$\be\!\cdot\!2<\al$}.  This is enough
  to prove the theorem by Fact~\ref{ords}.
  
  Let $T$ be an \llkt\ of order \be\ for some $K$.  If $(x_i)_1^n\in
  T$ let $F=F((x_i)_1^n)\subset S(X^{\ast})$ be a finite set which
  1-norms a $(1/2)$-net in $S([x_i]_1^n)$.  Choose by
  Lemma~\ref{lem-for-ix=wa} $T_{(x_i)_1^n,F}$ a block subtree of $T$
  of order \be\ contained in $X_F$.  Let $(a_k)$ be the collection of
  maximal nodes of $T$ and if $a_k=(x_i)_1^n$ let $T(k)=\{ a_k\cup x :
  x\in T_{(x_i)_1^n,F} \}$.  Thus $T'=T\cup(\cup_kT(k))$ is an
  $\ell_1$-$6K$-tree of order \mbox{$\be\!\cdot\!2$} as required, and
  hence \mbox{$\be\!\cdot\!2<\al$}.
\end{proof}

\begin{cor}\label{ibx=ix}
  Let $X$ have a basis.  If $I(X)\geq\w^{\w}$, then $I(X)=I_b(X)$.  
\end{cor}

\begin{proof}
  Let $\ap>\w$, and suppose $I(X)=\w^{\ap}$.  Then for every $\be$
  with $\w\leq\be<\al$ there exists $K$ such that
  \mbox{$I(X,K)\geq\w\!\cdot\!\w^{\bt}=\w^{\bt}$}.  Thus
  $\ib{X,K'}\geq\w^{\bt}$ for some $K'$ by
  Proposition~\ref{ind-to-block-ind}, and hence $\ib{X}>\w^{\bt}$ for
  every $\be<\al$.  If $\ap$ is a limit ordinal, then
  $\w^{\ap}=\sup_{\bt<\ap}\w^{\bt}$ and so $\ib{X}\geq\w^{\ap}=I(X)$.
  Otherwise $\al=\al'+1$, where $\al'\geq\w$ and $\ib{X}>\w^{\ap'}$,
  which implies $\ib{X}\geq\w^{\ap}=I(X)$ since $\ib{X}=\w^{\gm}$ for
  some $\gm$ by Theorem~\ref{ibx=wa}.  In either case we know that
  $I(X)\geq\ib{X}$ and so they are equal.
  
  If $I(X)=\w^{\w}$, then $I(X)>\w^{n+1}$ for every $ n\geq1$ and
  hence $\ib{X}>\w^n$ for every $ n\geq1$ by
  Proposition~\ref{ind-to-block-ind}.  Thus $\ib{X}\geq\w^{\w}=I(X)$
  and so $I(X)=\ib{X}$ as required.
\end{proof}

\begin{cor} \label{ix-ibx=1or0}
  If $I(X)=\w^n$, then $\ib{X}=\w^m$ where $m=n$ or $n-1$.
\end{cor}

\begin{proof}
  This follows from similar arguments to those for the previous
  corollary.
\end{proof}

\begin{rem}\label{which-ords-for-index}
  We collect together some notes about which values \a\ may take when
  $I(X)=\w^{\a}$.
\begin{enumerate}
\item If $X$ does not contain $\ell_1^n$'s uniformly, then
  $I(X)=\w=\ib{X}$.  Also, if $X$ contains $\ell_1^n$'s uniformly,
  then $I(X)\geq\w^2$.  It is easy to see that $\ib{c_0}=\w$ (where
  the block basis index is calculated with respect to the unit vector
  basis for $c_0$) and so $I(c_0)=\w^2$ by
  Corollary~\ref{ix-ibx=1or0}.  Thus the two ordinal indices may
  indeed differ.  In fact, by Remark~\ref{ind-of-sch-Ca} below, for
  each $n\geq1$ there exists a Banach space $X_n$ with
  $I_b(X_n)=I(X_n)=\omega^{n+1}$ and for each $n\geq1$ there exists a
  Banach space $Y_n$ with $I_b(Y_n)=\omega^n$ while
  $I(Y_n)=\omega^{n+1}$.
  
\item If $I(X)<\omega^{\omega}$, then it is possible for a space $X$
  to have two bases $(x_i)$ and $(y_i)$ with $I_b(X,(x_i))\neq
  I_b(X,(y_i))$.  Indeed, for each $n\geq1$ let $H_n$ be the span of
  the first $2^n$ Haar functions in $C(\Delta)$ (where $\Delta$ is the
  Cantor set on $[0,1]$); if $X=(\sum H_n)_{c_0}$, then $X\simeq c_0$.
  Thus, if $(x_i)$ is a basis for $X$ equivalent to the unit vector
  basis of $c_0$, then $I_b(X,(x_i))=\omega$.  However, if $(y_i)$ is
  the basis for $X$ consisting of the Haar bases for the $H_n$'s
  strung together, then, since each basis for $H_n$ admits a block
  basis of length $n$ which is $1$-equivalent to the unit vector basis
  of $\ell_1^n$, we obtain $I_b(X,(y_i))=\omega^2$.  By
  Corollary~\ref{ix-ibx=1or0} the block basis indices for different
  bases can only differ by a factor of \w.
  
\item We note here that there are some ordinals \al\ for which there
  are no spaces $X$ with index $I(X)=\w^{\ap}$.  In particular, if
  \al\ is a limit ordinal, then there is no space $X$ with
  $I(X)=\w^{\w^{\ap}}$.  Otherwise, let $I(X)=\w^{\w^{\ap}}$, then for
  all $\al'<\al$ there is some $K$ such that there exists an \llkt\ of
  order $\w^{\w^{\ap'}}$, which we may then refine to get an
  $\ell_1$-$(1+\ee)$-block subtree of order $\w^{\w^{\ap'}}$ for any
  $\ee>0$, by Remark~\ref{rems-after-thm}~(ii).  Hence $X$ contains a
  block basis tree of constant $2$ and order $\w^{\w^{\ap}}$ (taking
  the union of these trees) and so $I(X)\geq\w^{\w^{\ap}+1}$.
  
\item By Remark~\ref{ind-of-sch-Ca} below, for every $\alpha<\omega_1$
  there exists a Banach space $X$ with $I(X)=\omega^{\alpha+1}$, and
  by Theorem~\ref{ib-tsb=wbw} below, there exists a Banach space
  $Y=T({\mathcal S}_{\omega^{\alpha}},1/2)$ with
  $I(Y)=\omega^{\omega^{\alpha+1}}$.
  
\item If $X$ is asymptotic $\ell_1$ (see for example~\cite{otw} for
  the definition of this), then $\ib{X}\geq\w^{\w}$ and so
  $I(X)=\ib{X}$.

\end{enumerate}

\end{rem}

\begin{q}
  For which limit ordinals \al\ do there exist Banach spaces $X$ with
  index $I(X)=\w^{\ap}$?
\end{q}

We have already shown that there exist Banach spaces with index
$\omega^{\a+1}$ for every $\a<\omega_1$.  We have also shown that we
cannot have indices of the form $\w^{\w^{\ap}}$ for $\alpha$ a limit
ordinal, and that we do have spaces with index of the form
$\omega^{\omega^{\alpha+1}}$, but this leaves the question open for
all other limit ordinals.\vspace{\baselineskip}

This completes the first part of the section.  We now apply some of
these results and methods to calculating the $\ell_1$ index of some
Tsirelson spaces.

\begin{defn} \label{s-alpha} 
  Schreier sets of order \al, $\cs_{\ap}$ \cite{aa}.
\end{defn} Let
$E,F\subseteq{\mathbf N}, n\geq1$.  We write $E<F$ if $\max E<\min F$
and $n<E$ if $\{n\}<E$.  Let $\mathcal{M,N}$ be collections of finite
sets of integers and $K=(k_i)\subseteq{\mathbf N}$.  We define
\begin{multline*}
  \mathcal{M[N]}=\{\cup_1^kF_i : F_i\in {\mathcal N}\ (i=1,\dots,k) 
          \mbox{ and } \exists E=\{m_1,\dots,m_k\}\in {\mathcal M} \\
  \mbox{with } m_1\leq F_1<m_2\leq F_2<\dots<m_k\leq F_k;\ k\geq1\}
\end{multline*}
and ${\mathcal M}(K)=\{ \{ k_i : i\in E \} : E\in {\mathcal M} \}$.  

The Schreier sets, $\cs_{\ap}$ for each $\al<\w_1$ are defined
inductively as follows: Let $\cs_0=\{\{n\} :
n\geq1\}\cup\{\emptyset\},\ \cs_1=\{F\subset {\mathbf N} : |F|\leq
F\}=\cs_1[\cs_0]$.  If $\cs_{\ap}$ has been defined let
$\cs_{\ap+1}=\cs_1[\cs_{\ap}]$.  If \al\ is a limit ordinal with
$\cs_{\ap'}$ defined for each $\ap'<\ap$ choose an increasing sequence
of ordinals with $\ap=\sup_n\ap_n$ and let
$\cs_{\ap}=\cup_{n=1}^{\infty}\{ F\in\cs_{\ap_n} : n\leq F\}$.

For $n\geq1$ let $(\cs_{\ap})^n=\{F=\cup_1^nF_i : F_i\in\cs_{\ap},\ 
F_1<\dots<F_n\}$ and let $[\cs_{\ap}]^n=\cs_{\ap}[\dots[\cs_{\ap}]]$
($n$ times).  A sequence $(E_i)_1^n$ of finite subsets of integers is
$\cs_{\ap}$ admissible if $E_1<\dots<E_n$ and $(\min\,
E_i)_1^n\in\cs_{\ap}$.

Note that $(\cs_{\ap},\subseteq)$ forms a tree, ${\rm
  Tree}(\cs_{\ap})$, of order $\w^{\ap}$ and $[\cs_{\ap}]^n$ forms a
tree, ${\rm Tree}([\cs_{\ap}]^n)$, of order $\w^{\ap\cdot n}$ (see
eg.~\cite{aa}).

\begin{defn} \label{tsb} Tsirelson spaces, $\tsb$ \cite{a}.\end{defn}
We first define $c_{00}$ to be the linear space of all real sequences
with finite support, and let $(e_i)_1^{\infty}$ be the unit vector
basis of $c_{00}$.  If $E\subset{\mathbf N}$, then let $Ex=\sum_{i\in
  E}a_ie_i$.

Using the Schreier sets, Argyros defined the Tsirelson spaces, $T
(\cs_{\alpha},1/2)$, for $\alpha<\omega_1$.  He showed there exists a
norm \mbox{$\|\cdot\|$} on $c_{00}$ satisfying the implicit equation
\[ \|x\|=\max\left(\|x\|_{c_0},\frac{1}{2}
  \sup\Bigl\{\sum_{i=1}^n\|E_ix\| :   (E_i)_1^n\mbox{ is }
  \cs_{\a}\mbox{ admissible and }n\geq1\Bigr\}\right)\ . \]

The space $T(\cs_{\alpha},1/2)$ is the completion of
\mbox{$(c_{00},\|\cdot\|)$}.  The standard Tsirelson space $T$ (the
dual of Tsirelson's original space \cite{t}) is just $T(\cs_1,1/2)$
\cite{fj}.

\begin{defn} \label{Schreier-spaces}  
  Schreier spaces, $X_{\alpha}$.
\end{defn}
The Schreier spaces are generalizations of Schreier's example
\cite{sh}, first discussed in \cite{aa} and \cite{ao}.  They are
defined in a similar way to Tsirelson space; for each
$\alpha<\omega_1$ we define a norm on $c_{00}$ by:
\[ \left\|\sum a_ie_i\right\|_{\alpha}= \sup_{E\in{\mathcal 
    S}_{\alpha}}\Bigl|\sum_{i\in  E}a_i\Bigr|, \] 
and then the Schreier space $X_{\alpha}$ is the completion of
$(c_{00},\|\cdot\|_{\alpha})$.

\begin{thm} \label{ib-tsb=wbw}
  $\ib{T({\mathcal S}_{\a},1/2)}=\w^{\a\cdot\w}=I(T({\mathcal
    S}_{\a},1/2))$.
\end{thm}

\begin{prop} \label{prop-wbw}
  For each $\alpha<\omega_1$, for every $\varepsilon>0$, and for all
  $m\geq1$, there exists $n\geq1$ such that if $T$ is a block basis
  tree on a Banach space with a basis, and if ${\mathcal F}(T)=\{
  (\min(\spp x_i))_1^l : (x_i)_1^l\in T \}$ satisfies: $\forall F\in
  {\mathcal F}(T)\ \forall (a_i)_F\subset{\mathbf R}^{+}$ with
  $\sum_Fa_i=1$ there exists $G\in({\mathcal S}_{\a})^m$ such that
  $G\subset F$ and $\sum_Ga_i\geq\e$, then $o(T)\leq\w^{\a}\!\cdot\!
  n$.
\end{prop}
\begin{proof}
  We prove the result by induction on \a.  Let $\a=0$, pick $\ee>0,\ 
  m\geq1$ and choose $n$ so that $m/n<\ee$.  If $o(T)>\w^0\!\cdot\!
  n=n$, then there exists $ F\in {\mathcal F}(T),\ |F|>n$.  Now,
  setting $a_i=1/|F|$ for $i\in F$ gives $\sum_Fa_i=1$ but if
  $G\in({\mathcal S}_0)^m$, then $|G|=m$ and
  $\sum_Ga_i=m/|F|<m/n<\ee$, a contradiction.
  
  Suppose the result is true for \a.  To prove the case $\a+1$ first
  let $\ee>0$ be arbitrary and fix $m=1$.  Let $n>2/\ee$, and let $T$
  be a tree with \mbox{$o(T)\geq\w^{\a+1}\!\cdot\! n$}.  We may assume
  by Lemma~\ref{rep-in-min} that $T\simeq T(n,\w^{\a+1})$ and let
  $F:T\rightarrow T_n=\{a_1<\dots<a_n\}$ be the map
  $F^{-1}(a_1)=T_1(a_1)\mbox{ and
    }F^{-1}(a_i)=\cup_{n=1}^{\infty}T_n(a_i)\ (i>1)$ where
  $T_n(a_i)\simeq T_{\w^{\a+1}}$ and the $T_n(a_i)$'s are mutually
  incomparable.  Fix $m_1=1$ and $\ee_1<1/n$.  $T_1(a_1)$ has index
  $\w^{\a+1}>\w^{\a}\!\cdot\!k\ \forall k$ so $\exists F_1\in{\mathcal
    F}(T_1(a_1))\ \exists(a_i)_{F_1}\subset{\mathbf R}^{+}$ such that
  $\sum_{F_1}a_i=1$ and $\sum_Ga_i<\ee$ if $G\in({\mathcal
    S}_{\a})^{m_1}$.  Let $(x_i)_1^l\in T_1(a_1)$ be a node such that
  $F_1=(\min(\spp x_i))_1^l$.  Then there exists $i_2$ such that
  $T_{i_2}(a_2)>(x_i)_1^l$.
  
  Choose $m_2=\max(F_1)$ and $\ee_2<1/n$.  Repeating the process for
  the restricted tree $R(T_{i_2}(a_2))$ and $m_2,\ \ee_2$ up to
  $R(T_{i_n}(a_n))$ and $m_n,\ \ee_n$ we obtain $F_1<\dots<F_n,\ 
  (a_j)_{F_i}\subset{\mathbf R}^{+}$ such that $\sum_{F_i}a_j=1$ and
  $\sum_{G}a_j<\ee_i$ if $G\subseteq F_i$ and $G\in(\cs_{\a})^{m_i}$.
  Set $F=\cup_1^nF_i$ and $\overline{a}_j=\frac{1}{n}a_j$ for $j\in
  F$.  Let $G\subseteq F,\ G\in{\mathcal S}_{\a+1}$. Then
  $G=\cup_1^rG_j$ where $r\leq G_1<\dots<G_r$ and $G_j\in{\mathcal
    S_{\a}}$.  Let $i$ be least such that $G\cap F_i\neq\emptyset$
  then $r\leq\max(F_i)=m_{i+1}\leq m_l\ \forall l>i$. Hence if $l>i$,
  then $G\in(\cs_{\a})^{m_l}$ and so $\sum_{G\cap
    F_l}\overline{a}_j=\frac{1}{n}\sum_{G\cap F_l}a_j<\ee_l/n$;
  further $\sum_{G\cap
    F_i}\overline{a}_j\leq\sum_{F_i}\overline{a}_j=1/n$.  Thus
  \[ \sum_{G}\overline{a}_j<
  \frac{1}{n}(1+\ee_i+\dots+\ee_n)<\frac{2}{n}<\ee \] 
  as we had to show.
  
  For general $m>1$ we use the same construction, taking $n>2m/\ee$.
  Then, for $G\in(\cs_{\a+1})^m$ each set in ${\mathcal S}_{\a+1}$ can
  contribute at most $(1+\sum\ee_i)/n$ and hence we get the desired
  contradiction.
  
  Let $\a$ be a limit ordinal and suppose the result is true for each
  $\a'<\a$.  Let $(\a_i)$ be the increasing sequence of ordinals, with
  $\sup_i\a_i=\a$, which defines $\cs_{\a}$.  Let $\ee>0$, $m=1$, and
  choose $n>2/\e$.  Suppose \mbox{$o(T)\geq\w^{\a}\!\cdot\!n$}, and so
  assume $T\simeq T(n,\w^{\a})$; let $F:T\rightarrow
  T_n\equiv\{a_1<\dots<a_n\}$ be as before, but now with
  $T_n(a_i)\simeq T_{\w^{\a}}$.  From $\cf(T_1(a_1))$ select $F_1,\ 
  (a_i)_{F_1}\subset{\mathbf R}^{+}$ arbitrarily.  Let $(x_i)_1^l\in
  T_1(a_1)$ be a node such that $F_1=(\min(\spp x_i))_1^l$, then there
  exists $i\geq1$ such that $T_{i}(a_2)>(x_i)_1^l$; set
  $t_2=T_{i}(a_2)$.
  
  Now, the result is true for each $\a'<\a$, and
  \mbox{$o(R(t_2))>\w^{\a'}\!\cdot\!k$} for each $\a'<\a$ and every
  $k$, so there exists $F_2\in \cf(R(t_2))$, $m_2>\max F_1$, and
  $(a_j)_{F_2}$ such that $\sum_{F_2}a_j=1$ and every subset $G$ of
  $F_2$ which is also in $\cs_{\a_{m_2}}$ satisfies $\sum_Ga_j<\ee_2$,
  where $\e_2$ was chosen to be less than $1/n$.  Now, by \cite{otw},
  there exists $m$ such that if $G\geq m$ and $G\in\cs_{\a_i}$ for any
  $i<m_2$, then $G\in\cs_{\a_{m_2}}$.  Also, since $\w^{\a}$ is a
  limit ordinal, we may remove a finite number of the smallest nodes
  of $R(t_2)$ without changing the order of the tree and so we may
  choose $F_2\geq m$.
  
  We continue in this fashion, as before, to obtain $F_1<\dots<F_n,\ 
  (a_j)_{F_l}$ such that if $i\leq\max F_{l-1},\ G\in \cs_{\a_i},\ 
  G\subseteq F_l$, then $\sum_Ga_j<\ee_l<1/n$.  Set $F=\cup F_l$ and 
  $\overline{a}_j=\frac{1}{n}a_j$ for $j\in F$.  Let $G\in\cs_{\a}$,
  then there exists $j\geq1$ such that $G\in \cs_{\a_j}$ and $j\leq
  G$.  As before let $i$ be least such that $G\cap F_i\neq\emptyset$,
  then $j<m_l\ (l>i)$ and so $\sum_{G\cap
    F_l}\overline{a}_j=\frac{1}{n}\sum_{G\cap F_l}a_j<\ee_l/n$ and
  $\sum_{G\cap F_i}\overline{a}_j\leq\sum_{F_i}\overline{a}_j=1/n$.
  Thus
  \[ \sum_{G}\overline{a}_j<\frac{1}{n}(1+\ee_i+\dots+\ee_n)<
  \frac{2}{n}<\ee \]
  giving the required contradiction.
  
  The case for $m>1$ proceeds along similar lines as for the successor
  case; we just need to pick $n$ so that $m/n<\ee/2$.  This completes
  the proof of the proposition.
\end{proof}

\begin{proof}[Proof of Theorem~\ref{ib-tsb=wbw}]
  We first note that for each $n\geq1$, if $E\in[\cs_{\a}]^n$, then
  $\|\sum_{i\in E} a_ie_i\|\geq2^{-n}\sum|a_i|$, from the definition
  of the norm on $\tsb$, thus we may construct a block basis tree
  isomorphic to ${\rm Tree}([\cs_{\a}]^n)$.  As we noted in
  Definition~\ref{s-alpha} $o({\rm Tree}([\cs_{\a}]^n))=\w^{\a\cdot
    n}$, and hence $\ib{\tsb}\geq\w^{\a\cdot n}$ for each $n\geq1$,
  and so $\ib{\tsb}\geq\w^{\a\cdot\w}$.
  
  Now, suppose $\ib{\tsb}>\w^{\a\cdot\w}$, then there exists an
  \llkbbt\ $T$ of order $\w^{\a\cdot\w}$ and by Fact~\ref{ords}~(v) we
  may write $\w^{\a\cdot\w}=\w^{\w^{\theta+1}}$ for some
  $\theta<\w_1$.  This is one of the fixed points of our construction
  by Remark~\ref{rems-after-thm}~(ii).  Thus for every $\ee>0$ there
  exists an $\ell_1$-block subtree of $T$ with constant $1+\ee$ and
  order $\w^{\a\cdot\w}$, so we may assume $T$ has constant $1+\ee$
  where $\ee<1/10$.
  
  Let $m=1$ and choose $n$ from Proposition~\ref{prop-wbw}.  Since
  $o(T)>\w^{\a}\cdot n$ there exist $F\in{\mathcal F}(T),\ 
  F=\{n_1,\dots,n_l\}=(\min\,\spp x_i)_1^l$ for some $(x_i)_1^l\in T$
  and $(a_j)_F\subset{\mathbf R^+}$ such that $\sum_Fa_j=1$ and
  $\sum_Ga_j<\ee/3$ for each subset $G\subseteq F$ which is also in
  $\cs_{\a}$; set $x=\sum_{i=1}^la_{n_i}x_i$.  To calculate the norm
  of $x$ let $(E_i)_1^k$ be $\cs_{\a}$ admissible.  Let $I=\{i :
  \spp(x_i)\subseteq E_j \mbox{ for some }j \}$, let $J=\{ i\leq l :
  i\not\in I\text{ and }\spp(x_i)\cap E_j\neq\emptyset$ for some $j\}$
  and note that since $(E_i)_1^k$ is $\cs_{\a}$ admissible, there
  exist $A,B,C\in\cs_{\a}$ such that $\{ n_j : j\in J\}= A\cup B\cup
  C$.  Now
  \begin{align*}
    \frac{1}{2}\sum_{j=1}^k\|E_jx\| 
      & \leq \frac{1}{2}\sum_{i=1}^la_{n_i}\sum_{j=1}^k\|E_jx_i\|
        =\frac{1}{2}\left(\sum_{i\in I}a_{n_i}\sum_{j=1}^k\|E_jx_i\|+
        \sum_{i\in J}a_{n_i}\sum_1^k\|E_jx_i\|\right)\\
      & \leq  \frac{1}{2}\sum_{i\in I}a_{n_i}
        +\sum_{i\in J}a_{n_i}\frac{1}{2}\sum_{j=1}^k\|E_jx_i\|\\
      & \leq  \frac{1}{2}+\sum_{i\in J}a_{n_i}\|x_i\|\\
      & \leq  \frac{1}{2}+\sum_{j\in A\cup B\cup C}a_j\\
      & \leq  \frac{1}{2}+3\frac{\ee}{3}
  \end{align*}
  and hence $\|x\|\leq1/2+\ee$.  However $(x_i)_1^l\in T$, an
  $\ell_1$-$(1+\ee)$-tree and so $\|x\|\geq1/(1+\ee)$ a contradiction.
  Thus $\ib{\tsb}=\w^{\a\cdot\w}$.
\end{proof}

\begin{rem}\label{ind-of-sch-Ca}
  The authors have recently calculated the index of two other classes
  of Banach spaces.  In \cite{jo} it is shown that the index for
  $C(K)$, where $K$ is a countable compact metric space, is given by
  \[ I(C(\omega^{\omega^{\alpha}}))=\left\{
    \begin{array}{c} 
      \omega^{\alpha+2}\ (0\leq\alpha<\omega) \\
      \omega^{\alpha+1}\ (\w\leq\alpha)
    \end{array} 
  \right.  \]
  
  and for $1\leq\alpha<\omega_1$, and $X_{\alpha}$ the Schreier space
  for $\alpha$ (Definition~\ref{Schreier-spaces}), then
  $I(X_{\alpha})=\omega^{\alpha+1}$.
\end{rem}

\section{Final remarks}

As we noted in the introduction, Theorem~\ref{theorem} is false for
$1<p<\infty$.  This is a consequence of $\ell_p$ being arbitrarily
distortable \cite{os}.  In particular the following is true.

\begin{thm}\label{not-for-l-p}
  For each $p,\ 1<p<\infty$, and every $L\geq1$, there exist $K>1$ and
  $\al<\om_1$ such that for any $\be<\om_1$ there exists a Banach
  space $X$ which contains an $\ell_p$-$K$-tree on $X$ of order at
  least \be, but no $\ell_p$-$L$-tree of order \al.
\end{thm}

\begin{proof}
  Fix $L\geq1$; then since $\ell_p$ is arbitrarily distortable there
  exists a Banach space $X$ isomorphic to $\ell_p$ satisfying
  $d(Y,\ell_p)>2L$ for every subspace $Y$ of $X$.  Clearly, as $X$ is
  isomorphic to $\ell_p$, there exists some constant $K$ so that $X$
  contains an $\ell_p$-$K$-tree on $X$ of order \be\ for each
  $\be<\om_1$.  If the theorem is false, then for each $\al<\om_1$
  there would exist an $\ell_p$-$L$-tree on $X$ of order at least \al.
  This in turn would imply~\cite{b} that $X$ contains a subspace $Y$
  with $d(Y,\ell_p)\leq L$, contradicting our original assumption.
  This completes the proof.
\end{proof}

The finite version of Theorem~\ref{theorem} for $\ell_p$ is true, as
we mentioned in the introduction. From this and our construction of
$T_{\omega}$ (Definition~\ref{t-alpha-def}) it is easy to see that if
we have an $\ell_p$-$K $-tree $T$ of order $\omega$ on a Banach space
$X$, then there exists a block subtree of $T$ which is an
$\ell_p$-$(1+\varepsilon) $-tree of order $\omega$. Thus it seems
reasonable to ask the following question.

\begin{q}
  For which ordinals $\alpha$ is Theorem~\ref{theorem} true for
  $1<p<\infty$, and what is their supremum?
\end{q}

\begin{defn}\label{csm}
  $\ell_p$-${\mathcal S}_{\alpha}$-spreading models (${\mathcal
    S}_{\alpha}$-SMs)
\end{defn}

We extend the definition of the $\ell_1$-spreading models introduced
by Kiriakouli and Negrepontis~\cite{kn} to $\ell_p$ $(1\leq
p\leq\infty)$.  A sequence $(x_n)_{n=1}^{\infty}$ has an
$\ell_p$-${\mathcal S}_{\alpha}$-spreading model, for some $1\leq
p\leq\infty$, with constant $K$, if $(x_i)_{i\in
  F}\stackrel{K}{\sim}\mbox{uvb }\ell_p^{|F|}$ for every
$F\in{\mathcal S}_{\alpha}$, where ${\mathcal S}_{\alpha}$ is the
collection of Schreier sets of order \al\ introduced in
Section~\ref{index}.

We can refine the constant of an $\ell_1$-SM from $K$ to
$(1+\varepsilon)$ on a block basis as we did above for $\ell_1$-trees,
but the proof is much more straightforward.  We also note that these
spreading models are a stronger notion than $\ell_1$-trees.

We need the following result~\cite{otw}:

\begin{lem} \label{otw} 
  \emph{(OTW)} For each pair $\al,\be<\om_1$ there exists
  $N\subseteq{\mathbf N}$ such that ${\mathcal S}_{\alpha}[{\mathcal
    S}_{\beta}](N)\subseteq{\mathcal S}_{\beta+\alpha}$.
\end{lem}

\begin{thm} \label{refine-csms}
  For any $ K>1$, every $\varepsilon >0$, and each $\alpha<\omega_1$,
  there exists $\beta<\omega_1$ such that if $(x_n)$ is a normalized
  basic sequence having an $\ell_1$-${\mathcal S}_{\beta}$-SM with
  constant $K$, then there exists a normalized block basis $(y_n)$ of
  $(x_n)$ having an $\ell_1$-${\mathcal S}_{\alpha}$-SM with constant
  $1+\varepsilon$.
\end{thm}

\begin{proof} 
  This follows immediately from the following lemma.
\end{proof}

\begin{lem}  \label{lem:a2SM=>aSM}
  Let $(x_n)$ be a normalized basic sequence having an
  $\ell_1$-${\mathcal S}_{\ap\cdot2}$-SM with constant $K$. Then there
  exists a normalized block basis $(y_n)$ of $(x_n)$ having an
  $\ell_1$-${\mathcal S}_{\alpha}$-SM with constant $\scriptstyle
  \sqrt{K}$.
\end{lem}

\begin{proof}
  For fixed $\al<\om_1$ choose, by Lemma~\ref{otw},
  $N=(n_i)\subseteq{\mathbf N}$ such that ${\mathcal
    S}_{\alpha}[{\mathcal S}_{\alpha}](N)\subseteq{\mathcal
    S}_{\ap\cdot2}$ and consider the subsequence
  $(x_{n_i})_1^{\infty}$. We know that since ${\mathcal
    S}_{\ap\cdot2}(N)\subseteq{\mathcal S}_{\ap\cdot2}$,
  \[ \Bigl\|\sum_{i\in F}a_ix_{n_i}\Bigr\|\geq \frac{1}{K}\sum_{i\in
    F}|a_i|, \text{ for every } (a_i)\subset {\mathbf R}, \text{ and }
    F\in{\mathcal S}_{\ap\cdot2}\ .\] 
    
  If there exists $k\geq1$ such that
  \[  \Bigl\|\sum_{i\in E}a_ix_{n_i}\Bigr\|\geq
      \frac{1}{\sqrt{K}}\sum_{i\in E}|a_i|,   
  \text{ for every }
      (a_i)\subset {\mathbf R}, \text{ and each } E\in{\mathcal
      S}_{\alpha}\ \mbox{ with }E>k \]  
  then we are finished since $E\in {\mathcal S}_{\alpha}$ implies
  $E+k\in {\mathcal S}_{\alpha}\ (k\geq 1)$.
  
  Otherwise there exists a normalized block basis $(y_j)$ of
  $(x_{n_i})$ satisfying
  \[y_j=\sum_{i\in E_j}a_ix_{n_i},\ \sum_{i\in E_j}|a_i|>\sqrt{K}\]
  with $E_j\in {\mathcal S}_{\alpha}\mbox{ and }E_j<E_{j+1}$ for each
  $j\geq1$.  Now, for each $E\in\cs_{\alpha}$ the set $F=\cup_{j\in
    E}E_j$ is an element of ${\mathcal S}_{\alpha}[{\mathcal
    S}_{\alpha}](N)$, which in turn is contained in ${\mathcal
    S}_{\ap\cdot2}$. Thus we obtain
  \[ \Bigl\|\sum_Eb_jy_j\Bigr\|\geq\frac{1}{\sqrt{K}}\sum_{E}|b_j|,
  \text{ for every } 
  (b_j)\subset {\mathbf R}, \text{ and } E\in{\mathcal S}_{\alpha}\] 
  using James' argument as in the proof of Theorem~\ref{theorem}.
\end{proof}

\begin{rem}\label{last-rems}{\ }
  We note here some closing points for this section.
\begin{enumerate}
  
\item For every $\al<\omega_1$ there exists a Banach space
  $X_{\alpha}$ with an $\ell_1$-tree of order \al\ but $X_{\alpha}$
  has no $\ell_1$-spreading models.  In fact $X_{\alpha}$ can be taken
  to be reflexive with all normalized weakly null sequences having an
  $\ell_2$-$(1+\varepsilon)$ subsequence.

\begin{proof}
  We use a similar construction to Szlenk~\cite{s}.  Let
  $X_k=\ell_1^k\ (k\geq1)$.  If $\al<\omega_1$ is a limit ordinal and
  we have constructed $X_{\beta}$ for each $\be<\al$ let
  $X_{\alpha}=(\sum_{\beta<\alpha}X_{\beta})_{\ell_2}$.  Given
  $X_{\alpha}$ let $X_{\alpha+1}=(X_{\alpha}\oplus{\mathbf
    R})_{\ell_1}$.
\end{proof}

\item As for the $\ell_p$-trees, Theorem~\ref{refine-csms} is also
  true for $p=\infty$ and false for $1<p<\infty$.  This follows from
  the proof of Theorem~\ref{not-for-l-p}.
  
\item It follows from Lemma~\ref{lem:a2SM=>aSM} that if $(x_n)$ is a
  normalized basic sequence having an $\ell_1$-${\mathcal
    S}_{\w^{\a}}$-SM with any constant, then for every $\b<\w^{\a}$,
  and any $\e>0$, there exists a normalized block basis $(y_n)$ of
  $(x_n)$ having an $\ell_1$-${\mathcal S}_{\b}$-SM with constant
  $1+\varepsilon$.
\end{enumerate}
\end{rem}

\end{document}